\documentclass[reqno]{amsart}
\usepackage{setspace}
\usepackage{graphicx}
\usepackage{amsmath}
\usepackage{amsfonts}
\usepackage{amssymb}

\begin{document}

\centerline{\huge \bf Characterization of the Two-Dimensional}

\medskip
\centerline{\huge \bf  Five-Fold Translative Tiles}

\bigskip
\centerline{\bf Qi Yang and Chuanming Zong}

\vspace{1cm}
\centerline{\begin{minipage}{12.5cm}
{\bf Abstract.} In 1885, Fedorov discovered that a convex domain can form a lattice tiling of the Euclidean plane if and only if it is a parallelogram or a centrally symmetric hexagon. It is known that there is no other convex domain which can form a two-, three- or four-fold translative tiling in the Euclidean plane, but there are centrally symmetric convex octagons and decagons which can form five-fold translative tilings. This paper characterizes all the convex domains which can form five-fold translative tilings of the Euclidean plane, which consist of two classes of octagons and one class of decagons.
\end{minipage}}


\vspace{1cm}
\noindent
{\Large\bf 1. Introduction}

\bigskip
\noindent
In 1885, Fedorov \cite{fedo} proved that {\it a convex domain can form a lattice tiling in the plane if and only if it is a parallelogram or a centrally symmetric hexagon; a convex body can form a lattice tiling in the space if and only if it is a parallelotope, an hexagonal prism, a rhombic dodecahedron, an elongated dodecahedron, or a truncated octahedron}. As a generalized inverse problem of Fedorov's discovery, in 1900 Hilbert \cite{hilb} listed the following question in the second part of his 18th problem: {\it Whether polyhedra also exist which do not appear as fundamental regions of groups of motions, by means of which nevertheless by a suitable juxtaposition of congruent copies a complete filling up of all space is possible.} Try to verify Hilbert's problem in the plane, in 1917 Bieberbach suggested Reinhardt (see \cite{rein}) to determine all the two-dimensional convex tiles. However, to complete the list turns out to be challenging and dramatic. Over the years, the list has been successively extended by Reinhardt, Kershner, James, Rice, Stein, Mann, McLoud-Mann and Von Derau (see \cite{mann,zong14}), its completeness has been mistakenly announced several times! In 2017, M. Rao \cite{rao} announced a completeness proof based on computer checks.

Let $K$ be a convex body with (relative) interior ${\rm int}(K)$ and (relative) boundary $\partial (K)$, and let $X$ be a discrete set, both in $\mathbb{E}^n$. We call $K+X$ a {\it translative tiling} of $\mathbb{E}^n$ and call $K$ a {\it translative tile} if $K+X=\mathbb{E}^n$ and the translates ${\rm int}(K)+{\bf x}_i$ are pairwise disjoint. In other words, if $K+X$ is both a packing and a covering in $\mathbb{E}^n$. In particular, we call $K+\Lambda$ a {\it lattice tiling} of $\mathbb{E}^n$ and call $K$ a {\it lattice tile} if $\Lambda $ is an $n$-dimensional lattice. Apparently, a translative tile must be a convex polytope. Usually, a lattice tile is called a {\it parallelohedron}.

As one can predict that to determine the parallelohedra in higher dimensions is complicated. Through the works of Delone \cite{delo}, $\check{S}$togrin \cite{stog} and Engel \cite{enge}, we know that there are exact $52$ combinatorially different types of parallelohedra in $\mathbb{E}^4$. A computer classification for the five-dimensional parallelohedra was announced by Dutour Sikiri$\acute{\rm c}$, Garber, Sch$\ddot{\rm u}$rmann and Waldmann \cite{dgsw} only in 2015.

Let $\Lambda $ be an $n$-dimensional lattice. The {\it Dirichlet-Voronoi cell} of $\Lambda $ is defined by
$$C=\left\{ {\bf x}: {\bf x}\in \mathbb{E}^n,\ \| {\bf x}, {\bf o}\|\le \| {\bf x}, \Lambda \|\right\},$$
where $\| X, Y\|$ denotes the Euclidean distance between $X$ and $Y$. Clearly, $C+\Lambda $ is a lattice tiling and the Dirichlet-Voronoi cell $C$ is a parallelohedron. In 1908, Voronoi \cite{voro} made a conjecture that {\it every parallelohedron is a linear transformation image of the Dirichlet-Voronoi cell of a suitable lattice.} In $\mathbb{E}^2$, $\mathbb{E}^3$ and $\mathbb{E}^4$, this conjecture was confirmed by Delone \cite{delo} in 1929. In higher dimensions, it is still open.

To characterize the translative tiles is another fascinating problem. At the first glance, translative tilings should be more complicated than lattice tilings. However, the dramatic story had a happy end! It was shown by Minkowski \cite{mink} in 1897 that {\it every translative tile must be centrally symmetric}. In 1954, Venkov \cite{venk} proved that {\it every translative tile must be a lattice tile $($parallelohedron$)$} (see \cite{alek} for generalizations). Later, a new proof for this beautiful result was independently discovered by McMullen \cite{mcmu}.

Let $X$ be a discrete multiset in $\mathbb{E}^n$ and let $k$ be a positive integer. We call $K+X$ a {\it $k$-fold translative tiling} of $\mathbb{E}^n$ and call $K$ a {\it translative $k$-tile} if every point ${\bf x}\in \mathbb{E}^n$ belongs to at least $k$ translates of $K$ in $K+X$ and every point ${\bf x}\in \mathbb{E}^n$ belongs to at most $k$ translates of ${\rm int}(K)$ in ${\rm int}(K)+X$. In other words, $K+X$ is both a $k$-fold packing and a $k$-fold covering in $\mathbb{E}^n$. In particular, we call $K+\Lambda$ a {$k$-fold lattice tiling} of $\mathbb{E}^n$ and call $K$ a {\it lattice $k$-tile} if $\Lambda $ is an $n$-dimensional lattice. Apparently, a translative $k$-tile must be a convex polytope. In fact, similar to Minkowski's characterization, it was shown by Gravin, Robins and Shiryaev \cite{grs} that {\it a translative $k$-tile must be a centrally symmetric polytope with centrally symmetric facets.}

Multiple tilings was first investigated by Furtw\"angler \cite{furt} in 1936 as a generalization of Minkowski's conjecture on cube tilings. Let $C$ denote the $n$-dimensional unit cube. Furtw\"angler made a conjecture that {\it every $k$-fold lattice tiling $C+\Lambda$ has twin cubes. In other words, every multiple lattice tiling $C+\Lambda$ has two cubes sharing a whole facet.} In the same paper, he proved the two- and three-dimensional cases. Unfortunately, when $n\ge 4$, this beautiful conjecture was disproved by Haj\'os \cite{hajo} in 1941. In 1979, Robinson \cite{robi} determined all the integer pairs $\{ n,k\}$ for which Furtw\"angler's conjecture is false. We refer to Zong \cite{zong05,zong06} for an introduction account and a detailed account on this fascinating problem, respectively, to pages 82-84 of Gruber and Lekkerkerker \cite{grub} for some generalizations.

Let $P$ denote an $n$-dimensional centrally symmetric convex polytope, let $\tau (P)$ be the smallest integer $k$ such that $P$ can form a $k$-fold translative tiling in $\mathbb{E}^n$, and let $\tau^* (P)$ be the smallest integer $k$ such that $P$ can form a $k$-fold lattice tiling in $\mathbb{E}^n$. For convenience, we define $\tau (P)=\infty $ if $P$ can not form translative tiling of any multiplicity. Clearly, for every centrally symmetric convex polytope we have
$$\tau (P)\le \tau^*(P).$$

In 1994, Bolle \cite{boll} proved that {\it every centrally symmetric lattice polygon is a lattice multiple tile}. However, little is known about the multiplicity. Let $\Lambda $ denote the two-dimensional integer lattice, and let $P_8$ denote the octagon with vertices $(1,0)$, $(2,0)$, $(3,1)$, $(3,2)$, $(2,3)$, $(1,3)$, $(0,2)$ and $(0,1)$. As a particular example of Bolle's theorem, it was discovered by Gravin, Robins and Shiryaev \cite{grs} that {\it $P_8+\Lambda$ is a seven-fold lattice tiling of $\mathbb{E}^2$.}

In 2000, Kolountzakis \cite{kolo} proved that, if $D$ is a two-dimensional convex domain which is not a parallelogram and $D+X$ is a multiple tiling in $\mathbb{E}^2$, then $X$ must be a finite union of translated two-dimensional lattices. In 2013, a similar result in $\mathbb{E}^3$ was discovered by Gravin, Kolountzakis, Robins and Shiryaev \cite{gkrs}.

In 2017, Yang and Zong \cite{yz1,yz2} studied the multiplicity of the multiple translative tilings by proving the following results: {\it
Besides parallelograms and centrally symmetric hexagons, there is no other convex domain which can form a two-, three- or four-fold translative tiling in the Euclidean plane. However, there are particular octagons and decagons which can form five-fold translative tilings.} Meanwhile, Zong \cite{zong} characterized all the two-dimensional five-fold lattice tiles.

\medskip
This paper characterizes all the two-dimensional five-fold translative tiles by proving the following theorem.

\smallskip
\noindent {\bf Theorem 1.} {\it A convex domain can form a five-fold translative tiling of the Euclidean plane if and only if it is a parallelogram, a centrally symmetric hexagon, a centrally symmetric octagon {\rm (}under a suitable affine linear transformation{\rm )} with vertices
${\bf v}_1=\left( {3\over 2}-{{5\alpha }\over 4}, -2\right)$, ${\bf v}_2=\left( -{1\over 2}-{{5\alpha }\over 4}, -2\right)$, ${\bf v}_3=\left( {{\alpha }\over 4}-{3\over 2}, 0\right)$, ${\bf v}_4=\left( {{\alpha }\over 4}-{3\over 2}, 1\right)$, ${\bf v}_5=-{\bf v}_1$, ${\bf v}_6=-{\bf v}_2$, ${\bf v}_7=-{\bf v}_3$ and ${\bf v}_8=-{\bf v}_4$, where $0<\alpha <{2\over 3}$, or with vertices ${\bf v}_1=(2-\beta, -3),$ ${\bf v}_2=(-\beta, -3),$ ${\bf v}_3=(-2, -1),$ ${\bf v}_4=(-2, 1),$ ${\bf v}_5=-{\bf v}_1$, ${\bf v}_6=-{\bf v}_2$, ${\bf v}_7=-{\bf v}_3$ and ${\bf v}_8=-{\bf v}_4$, where $0<\beta \le 1$, or a centrally symmetric decagon with $\mathbf{u}_{1}=(0,1),$ $\mathbf{u}_{2}=(1,1),$ $\mathbf{u}_{3}=({3\over 2}, {1\over 2}),$ $\mathbf{u}_{4}=({3\over 2},0),$ $\mathbf{u}_{5}=(1,-{1\over 2}),$ $\mathbf{u}_{6}=-\mathbf{u}_{1},$ $\mathbf{u}_{7}=-\mathbf{u}_{2},$ $\mathbf{u}_{8}=-\mathbf{u}_{3},$ $\mathbf{u}_{9}=-\mathbf{u}_{4}$ and $\mathbf{u}_{10}=-\mathbf{u}_{5}$ as the middle points of its edges.}

\medskip\noindent
{\bf Remark 1.} It was shown by Zong \cite{zong} that a convex domain can form a five-fold lattice tiling of the Euclidean plane if and only if it is a parallelogram, a centrally symmetric hexagon, a centrally symmetric octagon (under a suitable affine linear transformation) with vertices ${\bf v}_1=(-\alpha , -{3\over 2})$, ${\bf v}_2=(1-\alpha , -{3\over 2})$, ${\bf v}_3=(1+\alpha , -{1\over 2})$, ${\bf v}_4=(1-\alpha , {1\over 2})$, ${\bf v}_5=-{\bf v}_1$, ${\bf v}_6=-{\bf v}_2$, ${\bf v}_7=-{\bf v}_3$ and ${\bf v}_8=-{\bf v}_4$, where $0<\alpha <{1\over 4},$ or with vertices
${\bf v}_1=(\beta , -2)$, ${\bf v}_2=(1+\beta , -2)$, ${\bf v}_3=(1-\beta , 0)$, ${\bf v}_4=(\beta , 1)$, ${\bf v}_5=-{\bf v}_1$, ${\bf v}_6=-{\bf v}_2$, ${\bf v}_7=-{\bf v}_3$, ${\bf v}_8=-{\bf v}_4$, where ${1\over 4}<\beta <{1\over 3}$, or a centrally symmetric decagon (under a suitable affine linear transformation) with ${\bf u}_1=(0, 1)$, ${\bf u}_2=(1, 1)$, ${\bf u}_3=({3\over 2}, {1\over 2})$, ${\bf u}_4=({3\over 2}, 0)$, ${\bf u}_5=( 1,-{1\over 2})$, ${\bf u}_6=-{\bf u}_1$, ${\bf u}_7=-{\bf u}_2$, ${\bf u}_8=-{\bf u}_3$, ${\bf u}_9=-{\bf u}_4$ and ${\bf u}_{10}=-{\bf u}_5$ as the middle points of its edges. In fact, all the two-dimensional convex five-fold translative tiles are five-fold lattice tiles. They take different representations just for the proof purpose.

\vspace{0.8cm}
\noindent
{\Large\bf 2. Preparation}

\bigskip\noindent
Let $P_{2m}$ denote a centrally symmetric convex $2m$-gon centered at the origin, let ${\bf v}_1$, ${\bf v}_2$, $\ldots$, ${\bf v}_{2m}$ be the $2m$ vertices of $P_{2m}$ enumerated in the clock order, and let $G_1$, $G_2$, $\ldots $, $G_{2m}$ be the $2m$ edges, where $G_i$ is ended by ${\bf v}_i$ and ${\bf v}_{i+1}$. For convenience, we write
$$V=\{{\bf v}_1, {\bf v}_2, \ldots, {\bf v}_{2m}\}$$
and
$$\Gamma=\{G_1, G_2, \ldots, G_{2m}\}.$$

Assume that $P_{2m}+X$ is a $\tau (P_{2m})$-fold translative tiling in $\mathbb{E}^2$, where $X=\{{\bf x}_1, {\bf x}_2, {\bf x}_3, \ldots \}$ is a discrete multiset with ${\bf x}_1={\bf o}$. Now, let us observe the local structure of $P_{2m}+X$ at the vertices ${\bf v}\in V+X$.

Let $X^{\bf v}$ denote the subset of $X$ consisting of all points ${\bf x}_i$ such that
$${\bf v}\in \partial (P_{2m})+{\bf x}_i.$$
Since $P_{2m}+X$ is a multiple tiling, the set $X^{\bf v}$ can be divided into disjoint subsets $X^{\bf v}_1$, $X^{\bf v}_2$, $\ldots ,$ $X^{\bf v}_t$ such that the translates in $P_{2m}+X^{\bf v}_j$ can be re-enumerated as $P_{2m}+{\bf x}^j_1$, $P_{2m}+{\bf x}^j_2$, $\ldots $, $P_{2m}+{\bf x}^j_{s_j}$ satisfying the following conditions:

\medskip
\noindent
{\bf 1.} {\it ${\bf v}\in \partial (P_{2m})+{\bf x}^j_i$ holds for all $i=1, 2, \ldots, s_j.$}

\smallskip\noindent
{\bf 2.} {\it Let $\angle^j_i$ denote the inner angle of $P_{2m}+{\bf x}^j_i$ at ${\bf v}$ with two half-line edges $L^j_{i,1}$ and $L^j_{i,2}$, where $L^j_{i,1}$, ${\bf x}^j_i-{\bf v}$ and $L^j_{i,2}$ are in clock order. Then, the inner angles join properly as
$$L^j_{i,2}=L^j_{i+1,1}$$
holds for all $i=1,$ $2,$ $\ldots ,$ $s_j$, where $L^j_{s_j+1,1}=L^j_{1,1}$.}

\medskip
For convenience, we call such a sequence $P_{2m}+{\bf x}^j_1$, $P_{2m}+{\bf x}^j_2$, $\ldots $, $P_{2m}+{\bf x}^j_{s_j}$ an {\it adjacent wheel} at ${\bf v}$. It is easy to see that
$$\sum_{i=1}^{s_j}\angle^j_i =2w_j\cdot \pi$$
hold for positive integers $w_j$. Then we define
$$\phi ({\bf v})=\sum_{j=1}^tw_j= {1\over {2\pi }}\sum_{j=1}^t\sum_{i=1}^{s_j}\angle^j_i$$
and
$$\varphi ({\bf v})=\sharp \left\{ {\bf x}_i:\ {\bf x}_i\in X,\ {\bf v}\in {\rm int}(P_{2m})+{\bf x}_i\right\}.$$

\medskip
Clearly, if $P_{2m}+X$ is a $\tau (P_{2m})$-fold translative tiling of $\mathbb{E}^2$, then
$$\tau (P_{2m})= \varphi ({\bf v})+\phi ({\bf v})\eqno (1)$$
holds for all ${\bf v}\in V+X$.

\smallskip
To prove Theorem 1, we need the following known results.

\bigskip
\noindent
{\bf Lemma 1 (Yang and Zong \cite{yz2}).} {\it Assume that $P_{2m}$ is a centrally symmetric convex $2m$-gon centered at the origin and $P_{2m}+X$ is a $\tau(P_{2m})$-fold translative tiling of the plane, where $m\geq4$. If $\mathbf{v}\in V+X$ is a vertex and $G\in \Gamma+X$ is an edge with $\mathbf{v}$ as one of its two ends, then there are at least $\lceil(m-3)/2\rceil$ different translates $P_{2m}+\mathbf{x}_{i}$ satisfying both
$$\mathbf{v}\in\partial(P_{2m})+\mathbf{x}_{i}$$
and}
$$G\backslash\{\mathbf{v}\}\subset {\rm int}(P_{2m})+\mathbf{x}_{i}.$$

\medskip\noindent
{\bf Lemma 2 (Yang and Zong \cite{yz2}).} {\it Assume that $P_{2m}$ is a centrally symmetric convex $2m$-gon centered at the origin, $P_{2m}+X$ is a translative multiple tiling of the plane, and $\mathbf{v}\in V+X$. Then we have
$$\phi ({\bf v})=\kappa\cdot {{m-1}\over 2}+\ell\cdot {1\over 2},$$
where $\kappa $ is a positive integer and $\ell$ is the number of the edges in $\Gamma+X$ which take ${\bf v}$ as an interior point.}

\medskip
\noindent
{\bf Lemma 3 (Yang and Zong \cite{yz2}).} {\it When $m\ge 6$, we have}
$$\tau (P_{2m})\ge 6.$$

\medskip\noindent
{\bf Lemma 4 (Bolle \cite{boll}).} {\it A convex polygon is a $k$-fold lattice tile for a lattice $\Lambda$ and some positive integer $k$ if and only if the following conditions are satisfied:

\noindent
{\bf 1.} It is centrally symmetric.

\noindent
{\bf 2.} When it is centered at the origin, in the relative interior of each edge $G$ there is a point of ${1\over 2}\Lambda $.

\noindent
{\bf 3.} If the midpoint of $G$ is not in ${1\over 2}\Lambda $ then $G$ is a lattice vector of $\Lambda $.}

\vspace{0.8cm}
\noindent {\Large \bf 3. Proof of Theorem 1}

\bigskip
\noindent{\bf Lemma 5.} {\it A centrally symmetric convex octagon $P_8$ is a five-fold translative tile if and only if it is, under an affine linear transformation, one with vertices ${\bf v}_1=\left( {3\over 2}-{{5\alpha }\over 4}, -2\right)$, ${\bf v}_2=\left( -{1\over 2}-{{5\alpha }\over 4}, -2\right)$, ${\bf v}_3=\left( {{\alpha }\over 4}-{3\over 2}, 0\right)$, ${\bf v}_4=\left( {{\alpha }\over 4}-{3\over 2}, 1\right)$, ${\bf v}_5=-{\bf v}_1$, ${\bf v}_6=-{\bf v}_2$, ${\bf v}_7=-{\bf v}_3$ and ${\bf v}_8=-{\bf v}_4$, where $0<\alpha <{2\over 3}$, or one with vertices
${\bf v}_1=(2-\beta, -3),$ ${\bf v}_2=(-\beta, -3),$ ${\bf v}_3=(-2, -1),$ ${\bf v}_4=(-2, 1),$ ${\bf v}_5=-{\bf v}_1$, ${\bf v}_6=-{\bf v}_2$, ${\bf v}_7=-{\bf v}_3$ and ${\bf v}_8=-{\bf v}_4$, where $0<\beta \le 1$.}

\bigskip
\noindent{\bf Proof.} Suppose that $X$ is a discrete subset of $\mathbb{E}^2$ and $P_{8}+X$ is a five-fold translative tiling of the plane. First of all, it follows from Lemma 1 that
$$\varphi ({\bf v})\ge \left\lceil {{4-3}\over 2}\right\rceil =1\eqno(2)$$
holds for all ${\bf v}\in V+X.$  On the other hand, by Lemma 2 we have
$$\phi ({\bf v})=\kappa \cdot {3\over 2}+\ell \cdot {1\over 2},\eqno(3)$$
where $\kappa $ is a positive integer and $\ell $ is a nonnegative integer. In fact $\ell$ is the number of the edges which take ${\bf v}$ as an interior point. Thus, to prove the lemma, it is sufficient to deal with the following four cases:

\medskip
\noindent
{\bf Case 1.} {\it $\phi(\mathbf{v})\geq 5$ holds for a  vertex $\mathbf{v}\in V+X$.} It follows by (1) and (2) that
$$\tau(P_{8})=\varphi(\mathbf{v})+\phi(\mathbf{v})\geq 6,\eqno(4)$$
which contradicts the assumption that $P_{8}+X$ is a five-fold translative tiling of the plane.

\medskip
\noindent
{\bf Case 2.} {\it $\phi(\mathbf{v})=4$ holds for a  vertex $\mathbf{v}\in V+X$.} It follows by (3) that $\ell\not= 0$ and therefore $\mathbf{v}\in {\rm int}(G)$ holds for some $G\in\Gamma+X$. Assume that ${\bf v}^*_1$ and ${\bf v}^*_2$ are the two ends of $G$. Applying Lemma 1 to $\{ {\bf v}^*_1, G\}$ and $\{ {\bf v}^*_2, G\}$, respectively, one can deduce that
$$\varphi(\mathbf{v})\geq 2$$
and therefore
$$\tau(P_{8})=\varphi(\mathbf{v})+\phi(\mathbf{v})\geq 6,\eqno(5)$$
which contradicts the assumption.

\medskip
\noindent
{\bf Case 3.} {\it $\phi(\mathbf{v})=3$ holds for a  vertex $\mathbf{v}\in V+X$.} Then (3) has and only has two groups of solutions
$\{\kappa , \ell\}=\{ 1, 3\}$ or $\{2, 0\}$.

\smallskip\noindent
{\bf Subcase 3.1.} $\{\kappa , \ell\}=\{ 1, 3\}$. Then, there are three edges $G'_1$, $G'_2$ and $G'_3$ in $\Gamma +X$ satisfying
$${\bf v}\in {\rm int}(G'_i),\quad i=1,\ 2,\ 3.$$
Next, we study the multiplicity by considering the relative positions of these edges.

\smallskip\noindent
{\bf Subcase 3.1.1.} $G'_{1}=G'_{2}=G'_{3}$. Assume that ${\bf v}^*_1$ and ${\bf v}^*_2$ are the two ends of $G'_1$. Then $X^{{\bf v}^*_1}$ has two identical points. By computing the angle sum of all the adjacent wheels at ${\bf v}^*_1$ it can be deduced that
$$\phi(\mathbf{v}_{1}^{*})\geq4.$$
Then, by Case 1 and Case 2 we get
$$\tau(P_{8})=\phi(\mathbf{v}_{1}^{*})+\varphi(\mathbf{v}_{1}^{*})\geq6,\eqno(6)$$
which contradicts the assumption.

\smallskip\noindent
{\bf Subcase 3.1.2.} {\it $G'_{2}=G'_{3}$ and $G'_{1}\nparallel G'_{2}$.} Then there are two adjacent wheels at $\mathbf{v}$, one
 has five translates $P_8+{\bf x}_1$, $P_8+{\bf x}_2$, $\ldots $, $P_8+{\bf x}_5$ and the other has two translates $P_8+{\bf x}'_1$ and $P_8+{\bf x}'_2$, as shown by Figure 1.

\begin{figure}[!ht]
\centering
\includegraphics[scale=0.53]{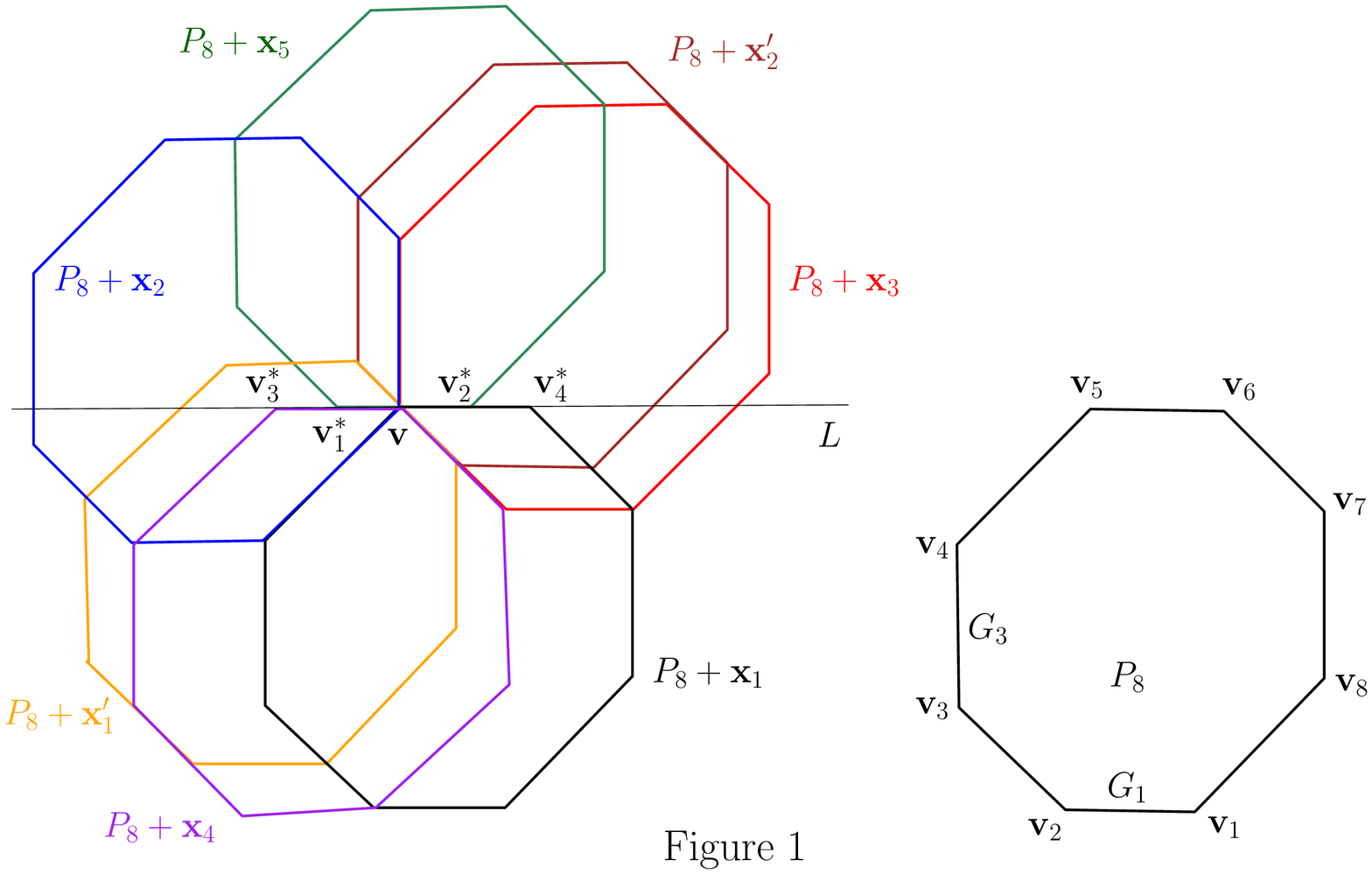}
\end{figure}

By re-enumeration we may assume that $\angle_{1},$ $\angle_{2},$ $\angle_{3}$ and $\angle_{4}$ are inner angles of $P_{8}$ and $\angle_{5}=\pi$, as shown by Figure 1. Guaranteed by linear transformation, we assume that the two edges $G_{1}$ and $G_{3}$ of $P_{8}$ are horizontal and vertical, respectively. Suppose that $G_{1}'=G_1+\mathbf{x}_{5}$. Let $\mathbf{v}_{1}^{*}$ and $\mathbf{v}_{2}^{*}$ be the two ends of $G_{1}'$, let $L$ denote the straight line determined by $\mathbf{v}_{1}^{*}$ and $\mathbf{v}_{2}^{*}$, let $G_{3}^{*}$ denote the edge of $P_{8}+\mathbf{x}_{4}$ lying on $L$ with two ends $\mathbf{v}$ and $\mathbf{v}_{3}^{*}$, and let $G_{4}^{*}$ denote the edge of $P_{8}+\mathbf{x}_{1}$ lying on $L$ with two ends $\mathbf{v}$ and $\mathbf{v}_{4}^{*}$.

By Lemma 1, there is a $\mathbf{x}^*_3\in X^{\mathbf{v}_{1}^{*}}$ such that $\mathbf{v}_{3}^{*}\in {\rm int}(P_{8})+\mathbf{x}^*_3$. Clearly, by the convexity of $P_8$, both $\mathbf{v}_{3}^{*}$ and $\mathbf{v}_{1}^{*}$ belong to ${\rm int}(P_{8})+\mathbf{x}'_1$. Thus, we have
$\mathbf{x}^*_3\neq \mathbf{x}'_1$. Meanwhile, since both $\mathbf{v}_{3}^{*}$ and $\mathbf{v}_{1}^{*}$ belong to ${\rm int}(P_{8})+\mathbf{x}_{2}$, we have $\mathbf{x}_2\neq \mathbf{x}^*_3$ and therefore
$$\varphi(\mathbf{v}_{3}^{*})\geq3.\eqno(7)$$
Then, the only chance to keep $P_8+X$ a five-fold tiling is $\phi(\mathbf{v}_{3}^{*})=2$. Similarly, one can deduce
$$\phi(\mathbf{v}_{1}^{*})=\phi(\mathbf{v}_{2}^{*})=\phi(\mathbf{v}_{3}^{*})=\phi(\mathbf{v}_{4}^{*})=2.\eqno(8)$$

By (3) it is easy to see that the local configuration of $P_8+X^{\bf v}$ is essentially unique when $\phi ({\bf v})=2$. In other words, it is determined by the one that ${\bf v}$ is not its vertex. Consequently, the set $X$ has four points ${\bf y}_1$, ${\bf y}_2$, ${\bf y}_3$ and ${\bf y}_4$ satisfying
$$\mathbf{v}_{1}^{*}=\mathbf{v}_{4}+\mathbf{y}_{1},\quad \mathbf{v}\in {\rm int}(P_{8})+\mathbf{y}_{1},\eqno(9)$$
$$\mathbf{v}_{2}^{*}=\mathbf{v}_{7}+\mathbf{y}_{2},\quad \mathbf{v}\in {\rm int}(P_{8})+\mathbf{y}_{2},\eqno(10)$$
$$\mathbf{v}_{3}^{*}=\mathbf{v}_{3}+\mathbf{y}_{3},\quad \mathbf{v}\in {\rm int}(P_{8})+\mathbf{y}_{3}\eqno(11)$$
and
$$\mathbf{v}_{4}^{*}=\mathbf{v}_{8}+\mathbf{y}_{4},\quad \mathbf{v}\in {\rm int}(P_{8})+\mathbf{y}_{4}.\eqno(12)$$

Clearly, by the convexity of $P_8$ we have $\mathbf{y}_{1}\neq\mathbf{ y}_{2}$, $\mathbf{y}_{1}\neq \mathbf{y}_{3}$ and $\mathbf{y}_{2}\neq \mathbf{y}_{4}$. For convenience, we write $\mathbf{v}_{i}=(x_{i},y_{i}).$ If $\mathbf{y}_{2}=\mathbf{y}_{3}$, then by (10) and (11) we have
$$y_{3}=y_{7}.\eqno(13)$$
If $\mathbf{y}_{1}=\mathbf{y}_{4}$, then by (9) and (12) we get
$$y_{4}=y_{8}.\eqno(14)$$
However, it is obvious that (13) and (14) can not hold simultaneously. Therefore, we still get
$$\varphi(\mathbf{v})\geq3$$
and therefore
$$\tau(P_{8})=\varphi(\mathbf{v})+\phi(\mathbf{v})\geq6,\eqno(15)$$
which contradicts the assumption.

\begin{figure}[!ht]
\centering
\includegraphics[scale=0.53]{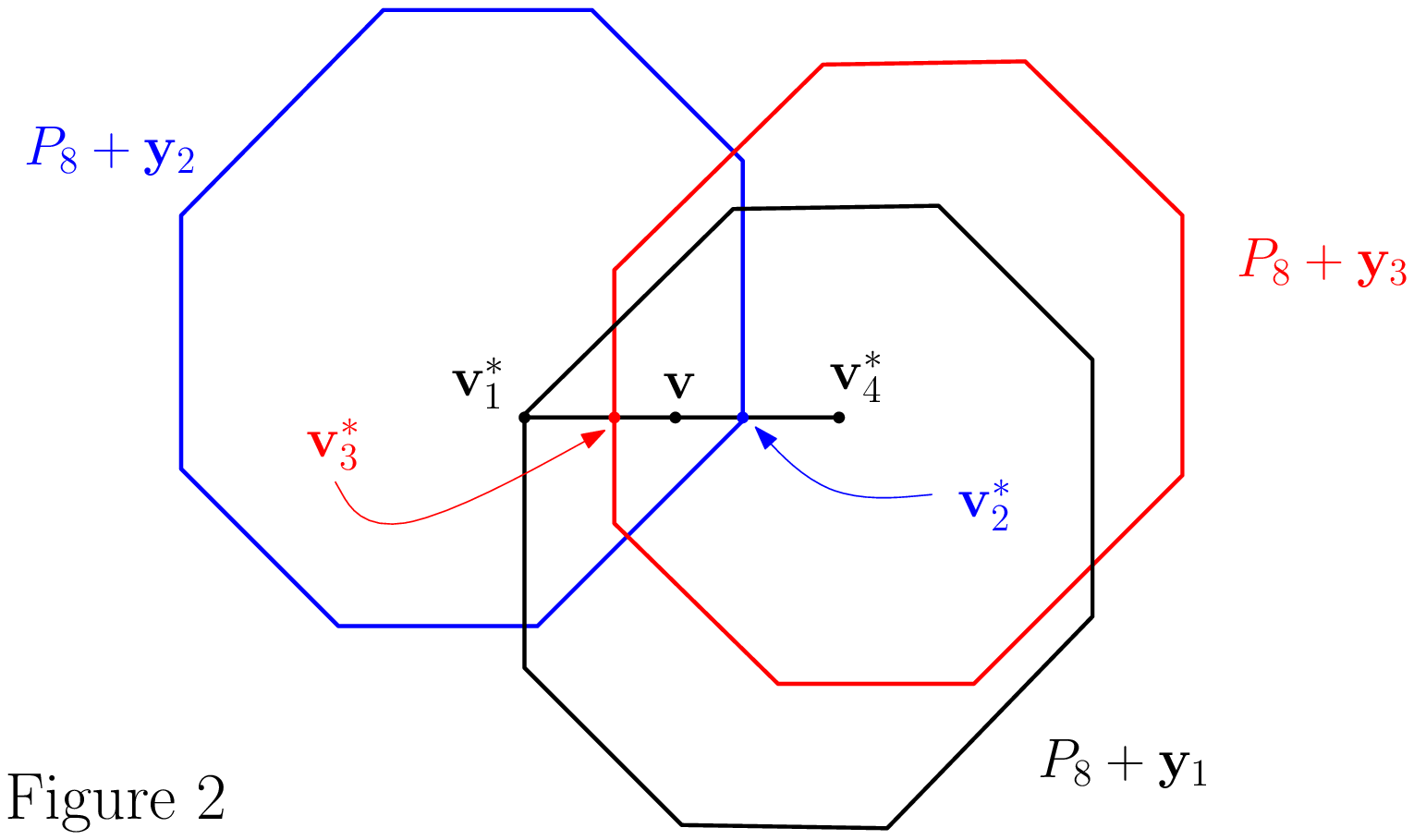}
\end{figure}

\smallskip
\noindent{\bf Subcase 3.1.3.} {\it $G'_1\neq G'_2$ and $G'_1\parallel G'_2$}. Let $\mathbf{v}_{1}^{*}$ and $\mathbf{v}_{2}^{*}$ be the two ends of $G'_{1}$, and let $\mathbf{v}_{3}^{*}$ and $\mathbf{v}_{4}^{*}$ be the two ends of $G'_{2}$. Without loss of generality, we suppose that $\mathbf{v}_{3}^{*}$ is between $\mathbf{v}_{1}^{*}$ and $\mathbf{v}_{2}^{*}$, as shown by Figure 2. By Lemma 1, $X$ has three points ${\bf y}_1$, ${\bf y}_2$ and ${\bf y}_3$ satisfying both
$${\bf v}^*_i\in \partial (P_8)+{\bf y}_i,\quad i=1,\ 2,\ 3$$
and
$${\bf v}\in {\rm int}(P_8)+{\bf y}_i,\quad i=1,\ 2,\ 3.$$
By the convexity of $P_8$ it is easy to see that these three points are pairwise distinct. Then, we get
$$\varphi(\mathbf{v})\geq3$$
and therefore
$$\tau(P_{8})=\phi(\mathbf{v})+\varphi(\mathbf{v})\geq6,\eqno(16)$$
which contradicts the assumption.

\smallskip
\noindent{\bf Subcase 3.1.4.} $G'_{1}\nparallel G'_{2},$ $G'_{1}\nparallel G'_{3}$ and $G'_{2}\nparallel G'_{3}$. By studying the angle sum at $\mathbf{v}$, it can be deduced that $P_{8}+X^\mathbf{v}$ is an adjacent wheel of seven translates.  Suppose that $\mathbf{x}_2\in X^\mathbf{v}$ and  $G'_{1}$ is an edge of $P_{8}+\mathbf{x}_2$. Since $G'_{1},$ $ G'_{2}$ and $G'_{3}$ are mutually non-collinear, $X^{\bf v}$ has two points
$\mathbf{x}_1$ and $\mathbf{x}_{3}$ such that $\mathbf{v}$ is a common vertex of both $P_{8}+\mathbf{x}_1$ and  $P_{8}+\mathbf{x}_3$, and $P_{8}+\mathbf{x}_2$ joins both $P_{8}+\mathbf{x}_1$ and $P_{8}+\mathbf{x}_3$ at non-singleton parts of $G'_{1}$, respectively. Let $\mathbf{v}^*_{1}$ and $\mathbf{v}^*_{2}$ be the two ends of $G'_{1}$, let $L$ denote the straight line determined by $\mathbf{v}^*_{1}$ and $\mathbf{v}^*_{2}$, let $G^*_{1}$ denote the edge of $P_{8}+\mathbf{x}_1$ lying on $L$ with ends $\mathbf{v}$ and $\mathbf{v}^*_{3}$, and let $G^*_{2}$ denote the edge of $P_{8}+\mathbf{x}_{3}$ lying on $L$ with ends $\mathbf{v}$ and $\mathbf{v}^*_{4}$, as shown in Figure 3.

\begin{figure}[!ht]
\centering
\includegraphics[scale=0.53]{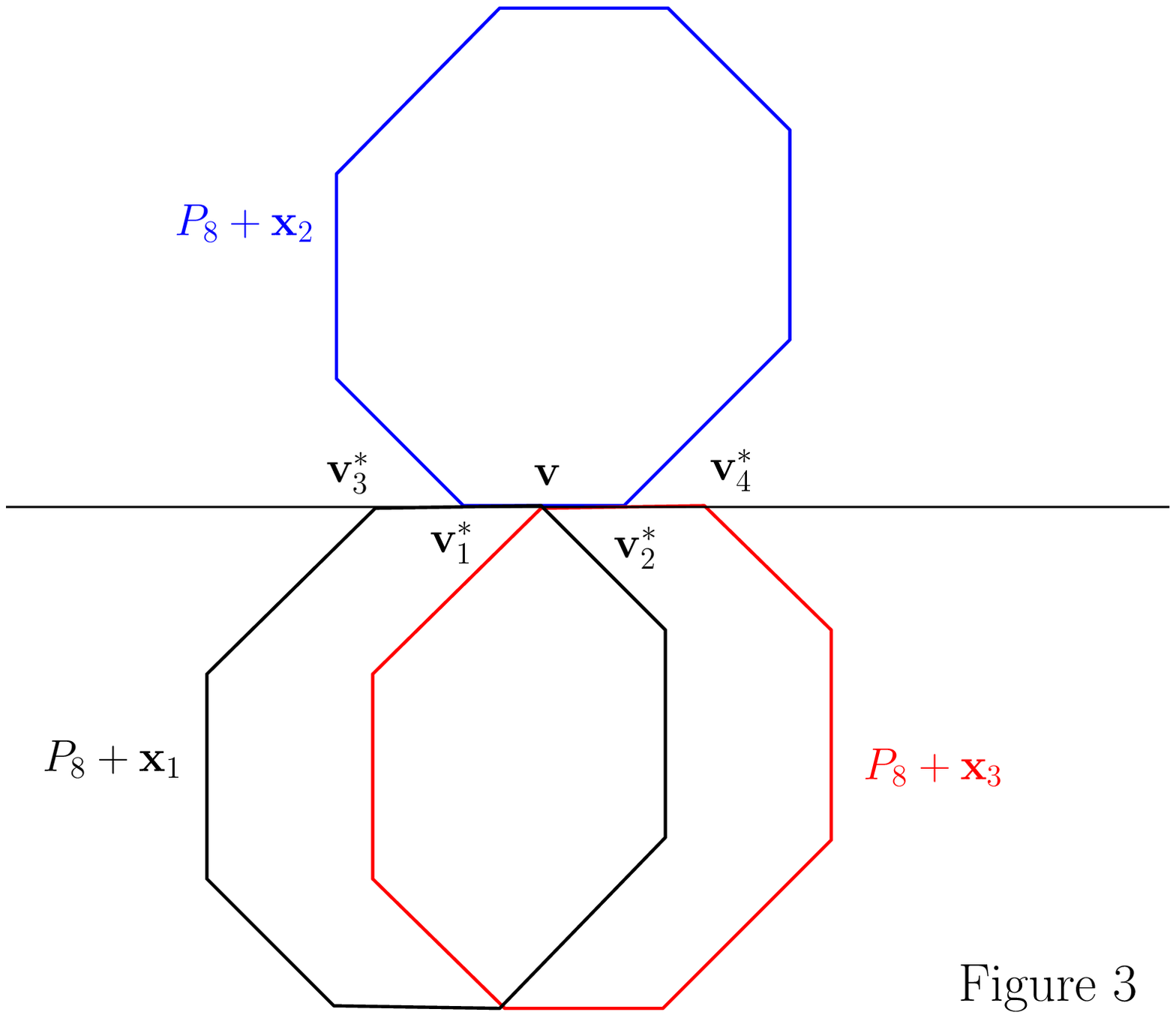}
\end{figure}

By studying the corresponding angles of the adjacent wheel at ${\bf v}$, it is easy to see that $P_{8}+X^\mathbf{v}$ has exact two translates which contain both ${\bf v}^*_1$ and ${\bf v}^*_3$ as interior points. On the other hand, by Lemma 1, $P_8+X^{{\bf v}^*_3}$ has at least one more
translate which contains ${\bf v}^*_1$ as an interior point. Thus, we have
$$\varphi ({\bf v}^*_1)\ge 3.$$
Then, the only chance to keep $P_8+X$ a five-fold tiling is $\phi(\mathbf{v}_{1}^{*})=2$. Similarly, one can deduce
$$\phi(\mathbf{v}_{1}^{*})=\phi(\mathbf{v}_{2}^{*})=\phi(\mathbf{v}_{3}^{*})=\phi(\mathbf{v}_{4}^{*})=2.\eqno(17)$$
By repeating the argument between (8) and (15), it can be deduced that
$$\tau(P_{8})=\varphi(\mathbf{v})+\phi(\mathbf{v})\geq6,\eqno(18)$$
which contradicts the assumption.

\medskip
\noindent
{\bf Subcase 3.2.} {\it $\{ \kappa, \ell\}=\{2,0\}$ holds at every vertex ${\bf v}\in V+X$.} Then $P_{8}+ X^{\mathbf{v}}$ is an adjacent wheel of eight translates $P_{8}+\mathbf{x}_{1},$ $P_{8}+\mathbf{x}_{2},$ $\ldots ,$ $P_{8}+\mathbf{x}_{8}$, as shown in Figure 4. Let $\mathbf{v}_{i}^{*}$ be the second vertex of $P_8+{\bf x}_i$ connecting to $\mathbf{v}$ by an edge. Since $\phi ({\bf v})=3$, every ${\bf v}^*_i$ is an interior point of exact two of these eight translates. Consequently, for every ${\bf v}^*_i$, there are two different translates $P_8+{\bf y}_i$
and $P_8+{\bf y}'_i$ in $P_8+X^{{\bf v}^*_i}$ both contain ${\bf v}$ as an interior point.

\begin{figure}[!ht]
\centering
\includegraphics[scale=0.53]{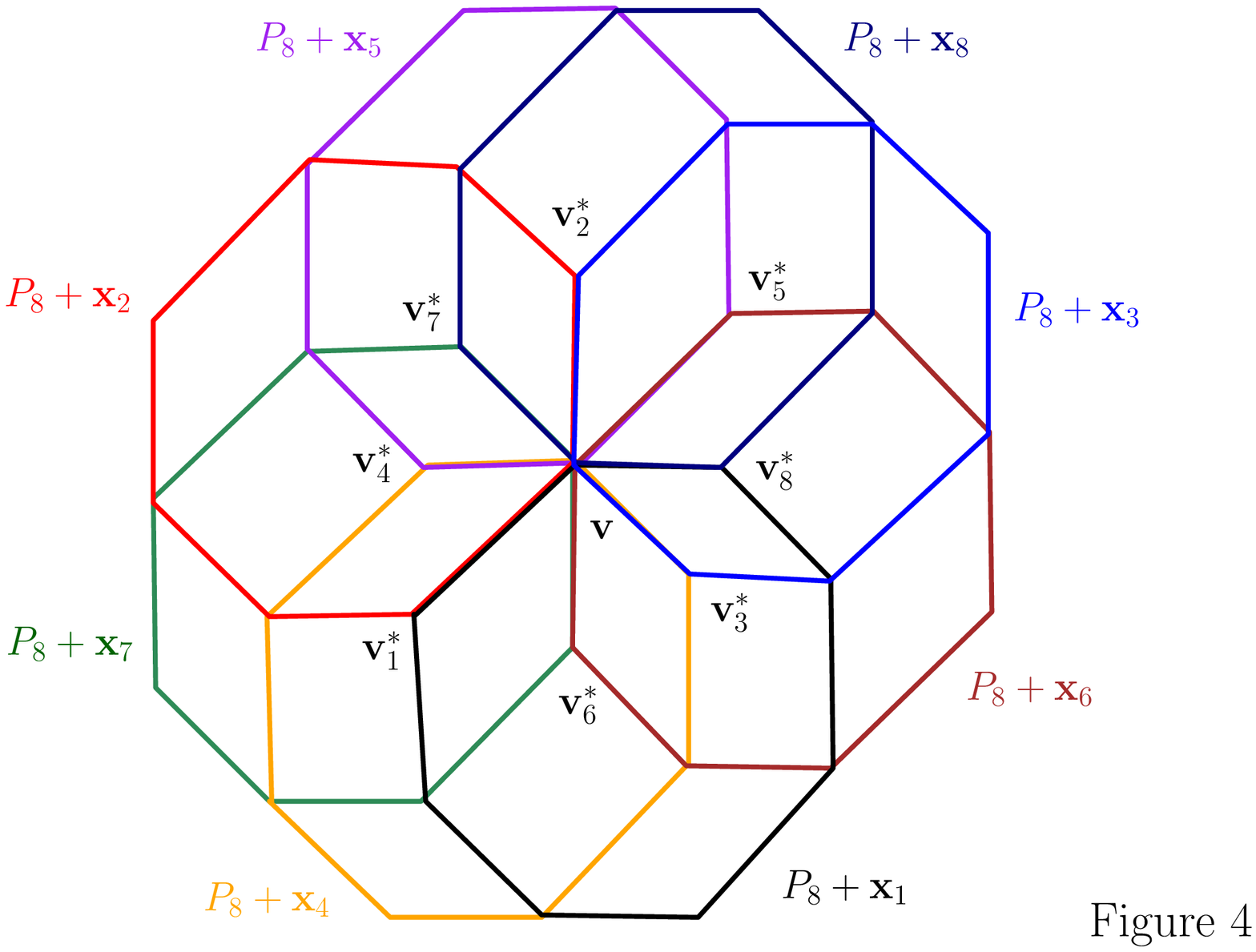}
\end{figure}

On the other hand, it can be easily deduced that there is only one point ${\bf x}\in X$ such that both ${\bf v}^*_1$ and ${\bf v}^*_2$ belong to
$\partial (P_8)+{\bf x}$ and ${\bf v}\in {\rm int}(P_8)+{\bf x}$. It is ${\bf v}^*_2-{\bf v}+{\bf x}_1$. Therefore, at least one of the two points ${\bf y}_2$ and ${\bf y}'_2$ is different from both ${\bf y}_1$ and ${\bf y}'_1$. Then, we get
$$\varphi(\mathbf{v})\geq3$$
and
$$\tau(P_{8})=\varphi(\mathbf{v})+\phi(\mathbf{v})\geq6,\eqno(19)$$
which contradicts the assumption.

\medskip
As a conclusion of the previous cases, if $P_8+X$ is a five-fold translative tiling, then $\phi(\mathbf{v})=2$ must hold at some $\mathbf{v}\in V+X$.

\medskip
\noindent
{\bf Case 4.} {\it $\phi(\mathbf{v})=2$ holds for a vertex $\mathbf{v}\in V+X$.} It follows by (3) that $\phi ({\bf v})=2$ holds if and only if $\kappa=1$ and $\ell =1$. In other words, $P_{8}+X^{\mathbf{v}}$ is an adjacent wheel of five translates. By re-enumeration we may assume that $\angle_{1},$ $\angle_{2},$ $\angle_{3}$ and $\angle_{4}$ are inner angles of $P_{8}$ and $\angle_{5}=\pi$, as shown  by Figure 5. Guaranteed by linear transformation, we assume that the edges $G_{1}$ and $G_{3}$ of $P_{8}$ are horizontal and vertical, respectively.

\begin{figure}[!ht]
\centering
\includegraphics[scale=0.53]{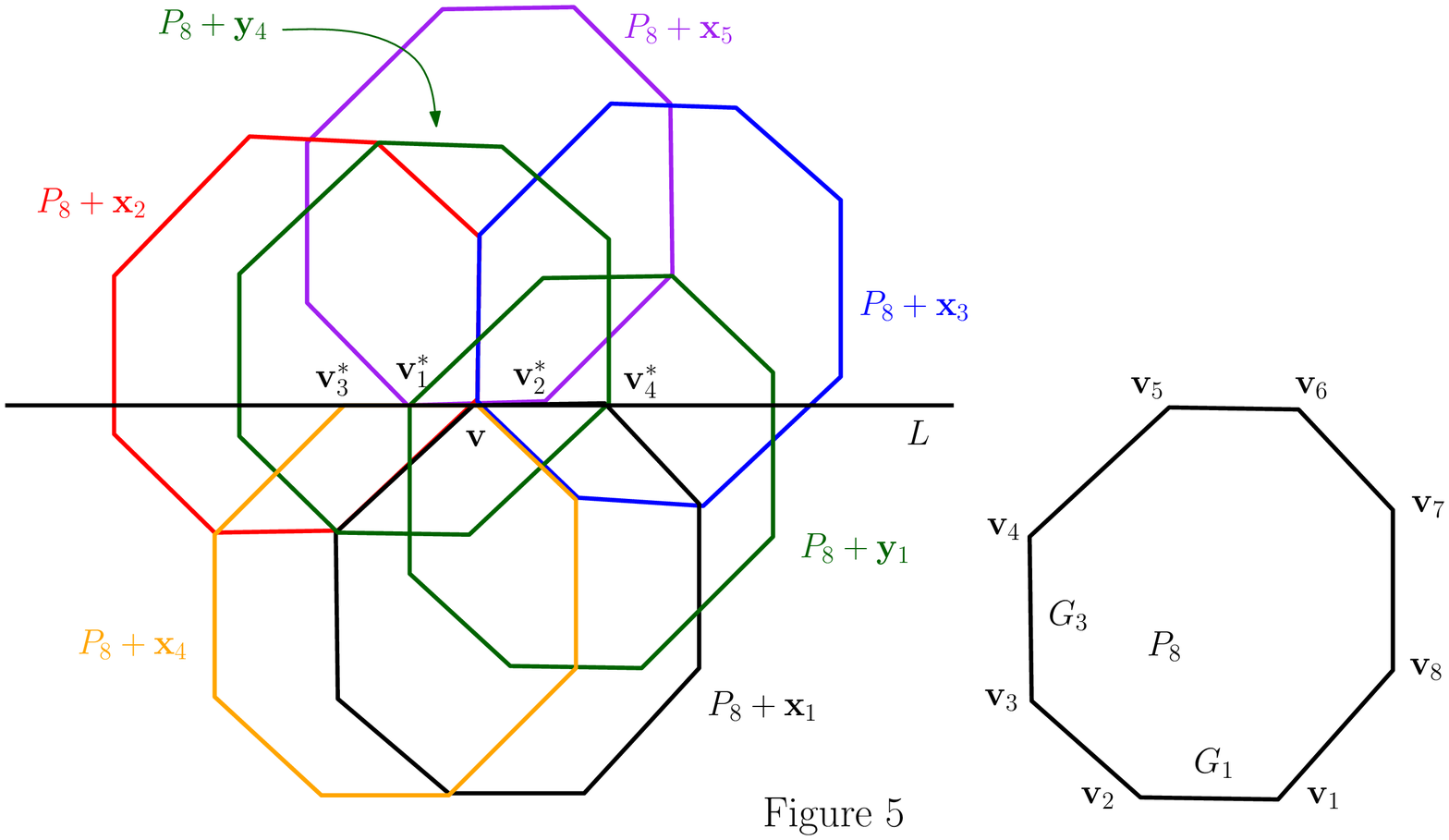}
\end{figure}

Let $G_{1}^{*}$ denote the edge of $P_{8}+\mathbf{x}_{5}$ such that $\mathbf{v}\in {\rm int}(G_{1}^{*})$ with two ends $\mathbf{v}_{1}^{*}$ and $\mathbf{v}_{2}^{*}$, let $L$ denote the straight line determined by $\mathbf{v}_{1}^{*}$ and $\mathbf{v}_{2}^{*}$, let $G_{3}^{*}$ denote the edge of $P_{8}+\mathbf{x}_{4}$ lying on $L$ with ends $\mathbf{v}$ and $\mathbf{v}_{3}^{*}$, and let $G_{4}^{*}$ denote the edge of  $P_{8}+\mathbf{x}_{1}$ lying on $L$ with ends $\mathbf{v}$ and $\mathbf{v}_{4}^{*}$. If  $\phi(\mathbf{v}_{1}^{*})=3$ or  $\phi(\mathbf{v}_{2}^{*})=3$, by Subcase 3.1 we have $\tau(P_{8})\geq6$, which contradicts the assumption that $P_{8}+X$ is a five-fold translative tiling. Thus we have $\phi(\mathbf{v}_{1}^{*})=2$ and  $\phi(\mathbf{v}_{2}^{*})=2$. Similarly, we have
$$\phi(\mathbf{v}_1^{*})=\phi(\mathbf{v}_2^{*})=\phi(\mathbf{v}_3^{*})=\phi(\mathbf{v}_4^{*})=2.\eqno(20)$$

Since the configuration of $P_8+X^{\bf v}$ is essentially unique if $\phi ({\bf v})=2$, by considering the wheel structures at ${\bf v}$, ${\bf v}^*_1$, ${\bf v}^*_2$,
${\bf v}^*_3$ and ${\bf v}^*_4$, there are four points ${\bf y}_1$, ${\bf y}_2$, ${\bf y}_3$ and ${\bf y}_4$ in $X$ satisfying
$$\mathbf{v}_{1}^{*}=\mathbf{v}_{4}+\mathbf{y}_{1},\quad \mathbf{v}\in {\rm int}(P_{8})+\mathbf{y}_{1},\eqno(21)$$
$$\mathbf{v}_{2}^{*}=\mathbf{v}_{7}+\mathbf{y}_{2},\quad \mathbf{v}\in {\rm int}(P_{8})+\mathbf{y}_{2},\eqno(22)$$
$$\mathbf{v}_{3}^{*}=\mathbf{v}_{3}+\mathbf{y}_{3},\quad \mathbf{v}\in {\rm int}(P_{8})+\mathbf{y}_{3}\eqno(23)$$
and
$$\mathbf{v}_{4}^{*}=\mathbf{v}_{8}+\mathbf{y}_{4},\quad \mathbf{v}\in {\rm int}(P_{8})+\mathbf{y}_{4}.\eqno(24)$$

By the convexity of $P_8$ it follows that $\mathbf{y}_{1}\neq \mathbf{y}_{2}$, $\mathbf{y}_{1}\neq \mathbf{y}_{3}$ and $\mathbf{y}_{2}\neq \mathbf{y_{4}}$. For convenience, we write $\mathbf{v}_{i}=(x_{i},y_{i}).$ If $\mathbf{y}_{1}=\mathbf{y}_{4}$, then by (21) and (24) we have
$$y_{4}=y_{8}.\eqno(25)$$
If $\mathbf{y}_{2}=\mathbf{y}_{3}$, then by (22) and (23) we have
$$y_{3}=y_{7}.\eqno(26)$$
It is obvious that (25) and (26) can not hold simultaneously. Therefore, we have either ${\bf y}_1\not= {\bf y}_4$ or ${\bf y}_2\not= {\bf y}_3$.

On the other hand, since $\phi ({\bf v})=2$, the three inequalities $\mathbf{y}_{3}\neq \mathbf{y}_{4}$, $\mathbf{y}_{2}\neq \mathbf{y}_{3}$ and $\mathbf{y}_{1}\neq \mathbf{y_{4}}$ can not hold simultaneously. Otherwise, it can be deduced that
$$\varphi({\bf v})\ge 4$$
and therefore
$$\tau(P_{8})=\varphi(\mathbf{v})+\phi(\mathbf{v})\geq 6,\eqno(27)$$
which contradicts the assumption. Since $\mathbf{y}_{1}=\mathbf{y}_{4}$ and $\mathbf{y}_{2}=\mathbf{y}_{3}$ are symmetric, it is sufficient to deal with two subcases.

\medskip
\noindent
{\bf Subcase 4.1.} $\mathbf{y}_{2}=\mathbf{y}_{3}$.  Let $\mathbf{v}'_1$ and $\mathbf{v}'_2$ be the two vertices of $P_{8}+\mathbf{x}_{2}$ that adjacent to $\mathbf{v}$, as shown in Figure 6.  By convexity it is easy to see that $\mathbf{v}'_1\in {\rm int}(P_{8})+\mathbf{y_{4}}$. Since $\mathbf{y}_{2}=\mathbf{y}_{3}$, we have $y_{3}=y_{7}$. Then $\mathbf{v}'_1$ is an interior point of $P_{8}+\mathbf{y_{2}}$ as well. Thus we get $\varphi(\mathbf{v}'_1)\geq 3$. If $\phi(\mathbf{v}'_1)\geq 3$, then we have
$$\tau (P_8)=\phi(\mathbf{v}'_1)+\varphi(\mathbf{v}'_1)\geq 6,\eqno(28)$$
which contradicts the assumption. Thus we must have $\phi(\mathbf{v}'_1)=2$. By lemma 1, there is a point $\mathbf{y_{5}}\in X^{\mathbf{v}'_1}$ such that $\mathbf{v}\in {\rm int}(P_{8})+\mathbf{y_{5}}$.

\begin{figure}[!ht]
\centering
\includegraphics[scale=0.53]{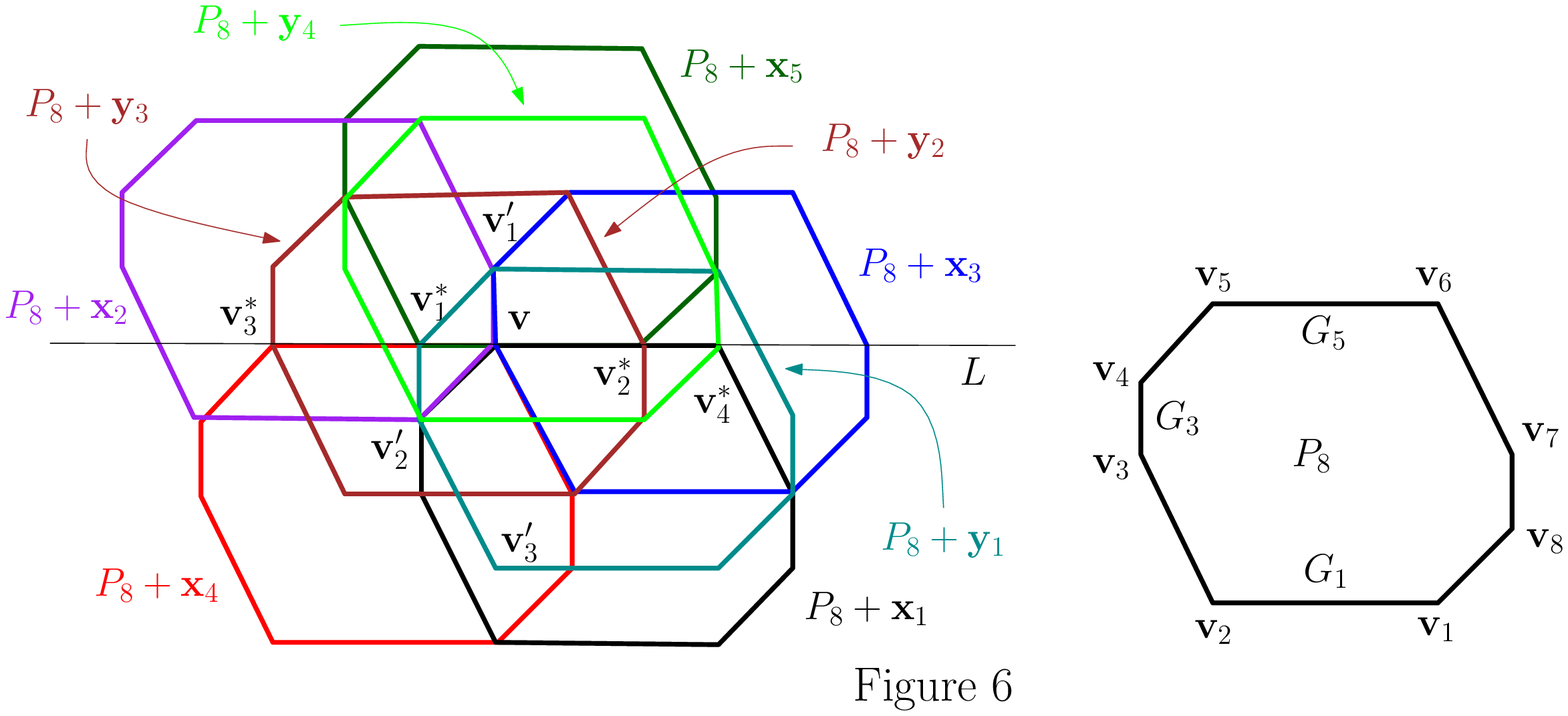}
\end{figure}

\noindent
{\bf Subcase 4.1.1.} {\it $\mathbf{v}'_1$ is a vertex of $P_{8}+\mathbf{y_{5}}$.} If $\mathbf{v}'_1$ is a vertex of $P_{8}+\mathbf{y_{1}}$, as shown by Figure 6, then one can deduce that $\mathbf{v}'_2$ is a vertex of $P_{8}+\mathbf{y_{1}}$. On the other hand, it follows by (24) that ${\bf v}'_2$ is a vertex of $P_8+{\bf y}_4$ as well. Therefore, we have $\phi ({\bf v}'_2)\not= 2$ and thus
$$\phi(\mathbf{v}'_2)\geq 3.\eqno(29)$$
Let $\mathbf{v}'_3$ denote the vertex $\mathbf{v}_{2}+\mathbf{y}_{1}$ of $P_{8}+\mathbf{y}_{1}$, as shown in Figure 6. By lemma 1, there is a point $\mathbf{z}\in X^{\mathbf{v}'_3}$ such that $\mathbf{v}'_2\in {\rm int}(P_{8})+\mathbf{z}$. Then it can be deduced that $\mathbf{v}'_3\in {\rm int}(P_{8})+\mathbf{x_{4}}$ and thus $\mathbf{z}\neq\mathbf{x}_{4}$. Since $y_{3}=y_{7}$, it can be shown that $\mathbf{v}'_3\notin P_{8}+\mathbf{y_{2}}$ and therefore $\mathbf{z}\neq\mathbf{y}_{2}$. In addition, we have
$$\mathbf{v}'_2\in \left( {\rm int}(P_{8})+\mathbf{x_{4}}\right)\cap \left( {\rm int}(P_{8})+\mathbf{y_{2}}\right).$$
Thus we have
$$\varphi(\mathbf{v}'_2)\geq3\eqno(30)$$
and consequently
$$\tau(P_{8})=\varphi(\mathbf{v}'_2)+\phi(\mathbf{v}'_2)\geq6,\eqno(31)$$
which contradicts the assumption.

If $\mathbf{v}'_1$ is not a vertex of $P_{8}+\mathbf{y_{1}}$, remembering the subcase assumption, then we have $\mathbf{y_{1}}\neq\mathbf{y_{5}}$. In fact, in this case all
${\bf y}_1$, ${\bf y}_3$, ${\bf y}_4$ and ${\bf y}_5$ are pairwise distinct. Thus we have
$$\varphi(\mathbf{v})\geq 4$$
and
$$\tau(P_{8})=\varphi(\mathbf{v})+\phi ({\bf v})\geq 6,\eqno(32)$$
which contradicts the assumption.

\medskip
\noindent
{\bf Subcase 4.1.2.} {\it $\mathbf{v}'_1$ is an interior point of an edge of $P_{8}+\mathbf{y_{5}}$.} It follows from the convexity of $P_8$ that $\mathbf{v}'_1 $ is an interior point of both $P_{8}+\mathbf{y}_{4}$ and $P_{8}+\mathbf{y}_{3}$. Therefore we have $\mathbf{y_{5}}\notin \{\mathbf{y}_3,\mathbf{y}_{4}\}$. If $\mathbf{y_{5}}\neq \mathbf{y}_1$, then all ${\bf y}_1$, ${\bf y}_3$, ${\bf y}_4$ and ${\bf y}_5$ are pairwise distinct. Thus we have
$$\varphi(\mathbf{v})\geq4$$
and
$$\tau(P_{8})=\varphi(\mathbf{v})+\phi ({\bf v})\geq6,\eqno(33)$$
which contradicts the assumption.

Thus, to avoid contradiction, we must have $\mathbf{y_{5}}=\mathbf{y}_{1}$. Notice that $\mathbf{v}'_1$ is an interior point of $P_{8}+\mathbf{x}_{5}$, and $P_{8}+\mathbf{y_{1}}$ has only two edges $G_{4}+\mathbf{y}_{1}$ and $G_{5}+\mathbf{y}_{1}$ which contain interior points of $P_{8}+\mathbf{x}_{5}$. Since $\phi(\mathbf{v}'_1)=2$, by studying the structure of the adjacent wheel at $\mathbf{v}'_1$, one can deduce that $\mathbf{v}'_1$ must be an interior point of $G_5+\mathbf{y}_{1}$. Then we have
$$y_5-y_{4}=y_{4}-y_{3}\eqno(34)$$
and
$$y_{3}-y_{2}=2 (y_{4}-y_{3}).\eqno(35)$$

Let $\mathbf{v}_{5}^{*}$ and $\mathbf{v}_{6}^{*}$ be the two ends of  $G_{5}+\mathbf{y}_1$. Suppose that $\mathbf{v}_{5}^{*}$ is on the left of $\mathbf{v}'_1$. By Lemma 1, there is a point $\mathbf{y}_{6}\in X^{\mathbf{v}_{5}^{*}}$ such that $\mathbf{v}'_1\in {\rm int}(P_{8}) +\mathbf{y}_{6}$.

It is obvious that $\mathbf{v}_{5}^{*}$ is an interior point of both $P_{8}+\mathbf{x_{5}}$ and $P_{8}+\mathbf{y_{2}}$. Thus we have $\mathbf{y}_{6}\notin\{\mathbf{y_{2}},\mathbf{x_{5}}\}$. If $\mathbf{v}_{5}^{*}$ is not lying on the boundary of $P_{8}+\mathbf{y_{4}}$, then we have $\mathbf{y_{4}}\neq\mathbf{y_{6}}$. Consequently, all ${\bf y}_2$, ${\bf y}_4$, ${\bf y}_6$ and ${\bf x}_5$ are pairwise distinct. Thus we have
$$\varphi(\mathbf{v}'_1)\geq4$$
and
$$\tau(P_{8})=\varphi(\mathbf{v}'_1)+\phi ({\bf v}'_1)\geq6,\eqno(36)$$
which contradicts the assumption.

To avoid the contradiction, the point $\mathbf{v}_{5}^{*}$ must belong to the boundary of $P_{8}+\mathbf{y_{4}}$. Furthermore, since the $y$-coordinate of $\mathbf{v}_{5}^{*}$ is equal to the $y$-coordinates of both $\mathbf{v}'_1$ and $\mathbf{v}_{3}+\mathbf{y}_{4}$, the point $\mathbf{v}_{5}^{*}$
must be the vertex $\mathbf{v}_{3}+\mathbf{y}_{4}$ of $P_{8}+\mathbf{y_{4}}$.

Let $v$ denote the $x$-coordinate of ${\bf v}$, and let $w_1$, $w_2$ and $w_3$ denote the $x$-coordinates of ${\bf v}_3+{\bf y}_4$, ${\bf v}^*_1$ and
${\bf v}^*_5$, respectively. First, by computing the $x$-coordinate of ${\bf v}^*_4$ in two ways we get
$$w_1+(x_7-x_6)+(x_6-x_5)+(x_5-x_4)=v+(x_6-x_5)$$
and thus
$$w_1 =v-(x_7-x_6)-(x_5-x_4).\eqno(37)$$
On the other hand, since ${\bf y}_2={\bf y}_3$, by computing the distance between ${\bf v}^*_3$ and ${\bf v}^*_4$ in two ways we get
$$(x_7-x_6)+(x_6-x_5)+(x_5-x_4)+v-w_2=2(x_6-x_5)$$
and thus
$$w_2=v+(x_7-x_6)-(x_6-x_5)+(x_5-x_4).$$
Since ${\bf v}^*_5$ is the left vertex of $G_5+{\bf y}_1$, we get
$$w_3=w_2+(x_5-x_4)=v+(x_7-x_6)-(x_6-x_5)+2(x_5-x_4).\eqno(38)$$
Then, ${\bf v}^*_5={\bf v}_3+{\bf y}_4$ implies $w_1=w_3$ and
$$2(x_7-x_6)+3(x_5-x_4)=x_6-x_5.\eqno(39)$$

In conclusion, recalling (34) and (35), a centrally symmetric octagon with $G_1$ horizontal, $G_3$ vertical and ${\bf y}_2={\bf y}_3$ is a five-fold translative tile only if
$$\left\{
\begin{array}{ll}
y_5-y_4=y_4-y_3,&\\
y_3-y_2=2(y_4-y_3),&\\
x_6-x_5=2(x_7-x_6)+3(x_5-x_4).&
\end{array}\right.\eqno(40)$$
Guaranteed by linear transformations, by choosing $y_4-y_3=1$, $x_6-x_5=2$ and $x_5-x_4=\alpha $ and keeping the symmetry in mind, one can deduce that the candidates are the octagons $D_8(\alpha)$ with vertices ${\bf v}_1=\left( {3\over 2}-{{5\alpha }\over 4}, -2\right)$, ${\bf v}_2=\left( -{1\over 2}-{{5\alpha }\over 4}, -2\right)$, ${\bf v}_3=\left( {{\alpha }\over 4}-{3\over 2}, 0\right)$, ${\bf v}_4=\left( {{\alpha }\over 4}-{3\over 2}, 1\right)$, ${\bf v}_5=-{\bf v}_1$, ${\bf v}_6=-{\bf v}_2$, ${\bf v}_7=-{\bf v}_3$ and ${\bf v}_8=-{\bf v}_4$, where $0<\alpha <{2\over 3}$.

Let $\Lambda (\alpha)$ denote the lattice generated by $\mathbf{u}_{1}=(2,0)$ and $\mathbf{u}_{2}=(1+{\alpha\over 2}, 1)$. It can be easily verified by Lemma 4 that $D_8(\alpha)+\Lambda (\alpha)$ is indeed a five-fold tiling of $\mathbb{E}^2$.

\medskip
\noindent
{\bf Subcase 4.2.} $\mathbf{y}_{3}=\mathbf{y}_{4}$. Then by (23) and (24) we have
$$y_{3}=y_{8}\eqno(41)$$
and
$$x_{8}-x_{3}=2(x_{1}-x_{2}).\eqno(42)$$

\begin{figure}[!ht]
\centering
\includegraphics[scale=0.53]{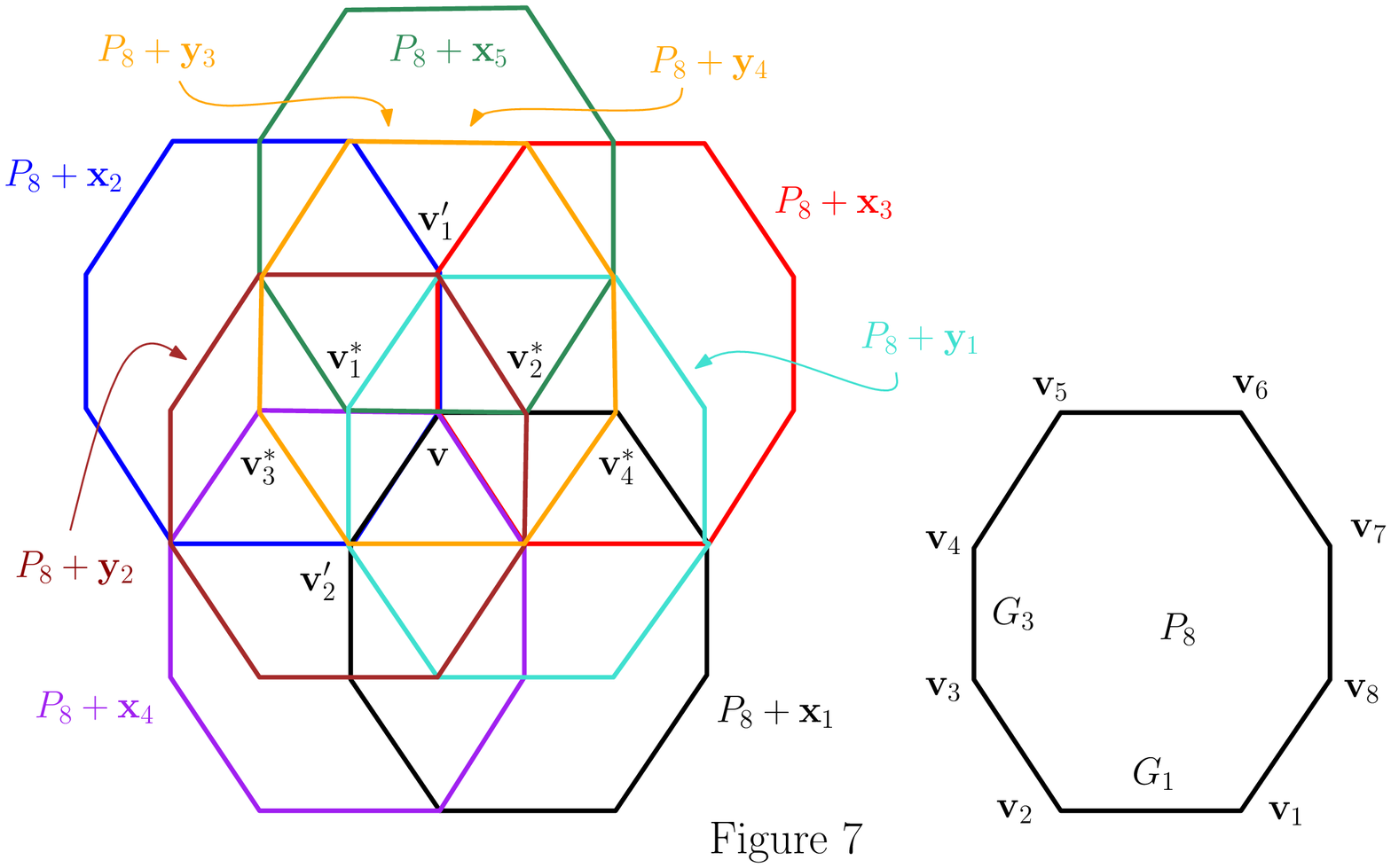}
\end{figure}

By Lemma 1, there is $\mathbf{y}_{5}\in X^{\mathbf{v}'_1}$ such that $\mathbf{v}\in {\rm int}(P_{8})+\mathbf{y}_{5}$. Since $y_{3}=y_{8}$, by convexity we have $\mathbf{v}'_1\in {\rm int}(P_{8})+\mathbf{y}_{3}$ and then $\mathbf{y}_{5}\neq \mathbf{y}_{3}$.

If $\phi(\mathbf{v}'_1)=3$ and $\mathbf{v}'_1\in {\rm int}(G)$ holds for some $G\in \Gamma+X$, by Subcase 3.1 we have
$$\tau(P_{8})=\varphi ({\bf v}'_1)+\phi ({\bf v}'_1)\geq6,\eqno(43)$$
which contradicts the assumption. Thus we have either $\phi(\mathbf{v}'_1)=2$ or $\phi(\mathbf{v}'_1)=3$ and $\mathbf{v}'_1$ is a vertex of $P_{8}+\mathbf{x}$ for all $\mathbf{x}\in X^{\mathbf{v}'_1}$.

If both $\mathbf{y}_{5}\neq\mathbf{y}_{1}$ and $\mathbf{y}_{5}\neq\mathbf{y}_{2}$ hold simultaneously, then we have
$$\varphi(\mathbf{v})\geq 4 $$
and therefore
$$\tau(P_{8})=\varphi ({\bf v})+\phi ({\bf v}) \geq6,\eqno(44)$$
which contradicts the assumption. So we must have either $\mathbf{y}_{5}=\mathbf{y}_{1}$ or $\mathbf{y}_{5}=\mathbf{y}_{2}$.

Suppose that $\mathbf{y}_{5}=\mathbf{y}_{1}$. If $\mathbf{v}'_1$ is a vertex of $P_{8}+\mathbf{y}_{1}$, then we have
$$y_{5}-y_{4}=y_{4}-y_{3}.\eqno(45)$$
If $\mathbf{v}'_1$ is an interior point of an edge of $P_{8}+\mathbf{y}_{1}$, then by Subcase 3.1 we have $\phi(\mathbf{v}'_1)=2$. By studying the structure of the adjacent wheel at $\mathbf{v}'_1$, one can deduce that $\mathbf{v}'_1$ must be an interior point of $G_{5}+\mathbf{y}_{1}$. Since $G_{5}+\mathbf{y}_{1}$ is horizontal, we also obtain (45).

In conclusion, recalling (41) and (42), a centrally symmetric octagon with $G_1$ horizontal, $G_3$ vertical and ${\bf y}_3={\bf y}_4$ is a five-fold translative tile only if
$$\left\{
\begin{array}{ll}
y_{3}=y_{8},&\\
y_{5}-y_{4}=y_{4}-y_{3},&\\
x_{8}-x_{3}=2(x_{1}-x_{2}).&
\end{array}\right.\eqno(46)$$
Guaranteed by linear transformation, by choosing $y_4-y_3=2$, $x_1-x_2=2$ and $x_6=\beta $ and keeping symmetry in mind, one can deduce that the candidates are the octagons $D_8(\beta)$ with vertices ${\bf v}_1=(2-\beta, -3),$ ${\bf v}_2=(-\beta, -3),$ ${\bf v}_3=(-2, -1),$ ${\bf v}_4=(-2, 1),$ ${\bf v}_5=-{\bf v}_1$, ${\bf v}_6=-{\bf v}_2$, ${\bf v}_7=-{\bf v}_3$ and ${\bf v}_8=-{\bf v}_4$, where $0<\beta \le 1$.

Let $\Lambda (\beta )$ denote the lattice generated by $\mathbf{u}_{1}=(2,0)$ and $\mathbf{u}_{2}=(1+{\beta \over 2}, 1)$. 
It can be easily verified by Lemma 4 that $D_8(\beta)+\Lambda (\beta)$ is indeed a five-fold tiling of $\mathbb{E}^2$.

Lemma 5 is  proved. \hfill{$\Box$}

\bigskip
\noindent
{\bf Lemma 6.} Let $P_{10}$ be a centrally symmetric decagon centered at the origin and let $X$ be a discrete multiset of $\mathbb{E}^2$. If $P_{10}+X$ is a five-fold translative tiling of $\mathbb{E}^2$, then it is  a five-fold lattice tiling of $\mathbb{E}^2$.

\medskip
\noindent {\bf Proof.} Let ${\bf v}_1$, ${\bf v}_2$, $\ldots $, ${\bf v}_{10}$ denote the ten vertices of $P_{10}$ enumerated in the clock order, let $G_i$ denote the edge with ends ${\bf v}_i$ and ${\bf v}_{i+1}$, and let ${\bf u}_i$ denote the middle point of $G_i$. Suppose that $X$ is a discrete subset of $\mathbb{E}^2$ and $P_{10}+X$ is a five-fold translative tiling of the plane. First of all, it follows from Lemma 1 that
$$\varphi ({\bf v})\ge \left\lceil {{5-3}\over 2}\right\rceil =1\eqno(47)$$
holds for every ${\bf v}\in V+X.$  On the other hand, by Lemma 2 we have
$$\phi ({\bf v})=\kappa \cdot 2+\ell \cdot {1\over 2},\eqno(48)$$
where $\kappa $ is a positive integer and $\ell $ is the number of the edges which contain ${\bf v}$ as a relative interior point.

Now we proceed to show that $\phi(\mathbf{v})=2$ holds for every vertex $\mathbf{v}\in V+X$ by dealing with the following two cases.

\medskip
\noindent{\bf Case 1.} {\it $\ell =0$ holds for all vertices $\mathbf{v}\in V+X$.} Then it follows by (48) that $\phi(\mathbf{v})$ can take only two values, two or four.

If $\phi(\mathbf{v})=4$, the local arrangements $P_{10}+X^{\bf v}$ can be divided into two adjacent wheels, each contains five translates.
Suppose that $P_{10}+{\bf x}_1$, $P_{10}+{\bf x}_2$, $\ldots$, $P_{10}+{\bf x}_5$ is such a wheel at ${\bf v}$ and $\mathbf{v}=\mathbf{v}_{k}+\mathbf{x}_{1}$. Then, the wheel can be determined by $P_{10}+\mathbf{x}_{1}$ explicitly as following:
$$\mathbf{v}=\mathbf{v}_{k+4}+\mathbf{x}_{2},\quad G_{k+4}+{\bf x}_2=G_{k-1}+\mathbf{x}_{1},$$
$$\mathbf{v}=\mathbf{v}_{k+8}+\mathbf{x}_{3},\quad G_{k+8}+{\bf x}_3=G_{k+3}+\mathbf{x}_{2},$$
$$\mathbf{v}=\mathbf{v}_{k+2}+\mathbf{x}_{4},\quad G_{k+2}+{\bf x}_4=G_{k+7}+\mathbf{x}_{3},$$
$$\mathbf{v}=\mathbf{v}_{k+6}+\mathbf{x}_{5},\quad G_{k+6}+{\bf x}_5=G_{k+1}+\mathbf{x}_{4},$$
where $\mathbf{v}_{10+i}=\mathbf{v}_{i}$ and $G_{10+i}=G_i$.

\begin{figure}[!ht]
\centering
\includegraphics[scale=0.5]{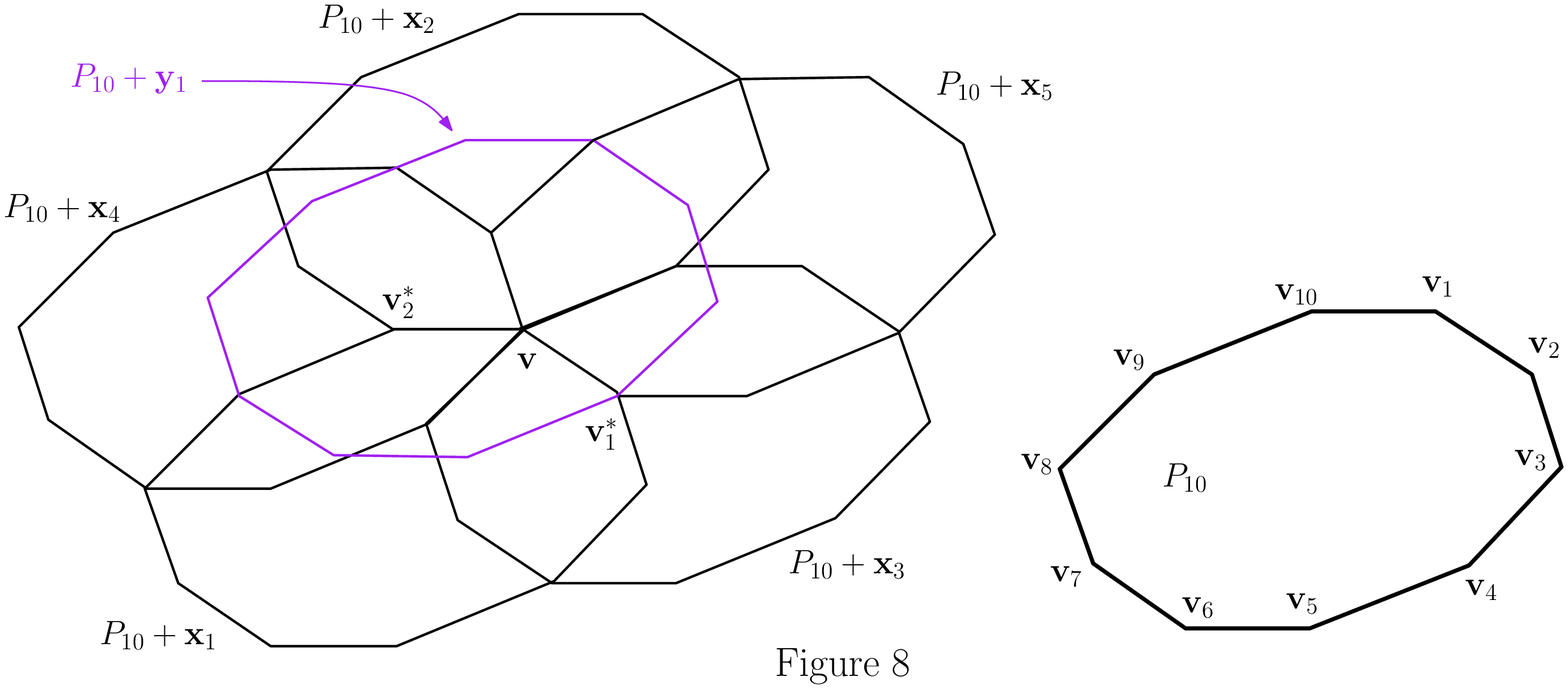}
\end{figure}

Without loss of generality, as shown by Figure 8, we take $\mathbf{v}=\mathbf{v}_{1}+\mathbf{x}_{1},$ $\mathbf{v}_{1}^{*}=\mathbf{v}_{2}+\mathbf{x}_{1}$ and $ \mathbf{v}_{2}^{*}=\mathbf{v}_{10}+\mathbf{x}_{1}$. By Lemma 1, for each ${\bf v}^*_i$ there is a ${\bf y}_i\in X^{{\bf v}^*_i}$ such that
$${\bf v}\in {\rm int}(P_{10})+{\bf y}_i.$$
In fact, by the previous analysis, we have ${\bf y}_1={\bf v}^*_1-{\bf v}_4$ and therefore by convexity
$${\bf v}^*_2\in {\rm int}(P_{10})+{\bf y}_1.$$
Thus, the two points ${\bf y}_1$ and ${\bf y}_2$ are different. Then we have
$$\varphi ({\bf v})\ge 2$$
and
$$\tau (P_{10})=\varphi ({\bf v})+\phi ({\bf v})\ge 6,\eqno(49)$$
which contradicts the assumption.

This means that, in this case $\phi(\mathbf{v})=2$ must hold for all ${\bf v}\in V+X.$

\medskip\noindent
{\bf Case 2.} {\it $\ell \not=0$ holds at a vertex $\mathbf{v}\in V+X$.} In other words, there is an edge $G\in \Gamma +X$ such that ${\bf v}\in {\rm int}(G)$. Clearly, by (48) we have $\phi ({\bf v})\ge 3.$

Suppose that $\mathbf{v}_{1}^{*}$ and $\mathbf{v}_{2}^{*}$ are the two ends of $G$. By Lemma 1, there are two different points $\mathbf{x}_{1}\in X^{\mathbf{v}_{1}^{*}}$ and $\mathbf{x}_{2}\in X^{\mathbf{v}_{2}^{*}}$ such that
$$\mathbf{v}\in \left({\rm int}(P_{10})+\mathbf{x}_{1}\right)\cap \left({\rm int}(P_{10})+\mathbf{x}_{2}\right).$$
Then we have $\varphi(\mathbf{v})\geq2$. If $\phi ({\bf v})\ge 4$, one can deduce that
$$\tau (P_{10})=\varphi(\mathbf{v})+\phi(\mathbf{v})\geq6,\eqno(50)$$
which contradicts the assumption.

If $\phi(\mathbf{v})=3$, by (48) one can deduce that $P_{10}+X^{\mathbf{v}}$ consists of seven translates $P_{10}+\mathbf{x}_{1},$ $P_{10}+\mathbf{x}_{2},$ $\ldots,$ $P_{10}+\mathbf{x}_{7}$, and there is another $G'\in \Gamma+X$ which contains $\mathbf{v}$ as an interior point. Suppose that $G$ is an edge of $P_{10}+\mathbf{x}_{6}$ and $G'$ has two ends $\mathbf{v}_{5}^{*}$ and $\mathbf{v}_{6}^{*}$. We deal with three subcases.

\smallskip\noindent
{\bf Subcase 2.1.} {\it $G'||G$ and $G'\neq G$.} Without loss of generality, we assume that $\mathbf{v}_{5}^{*}$ is between $\mathbf{v}_{1}^{*}$ and $\mathbf{v}_{2}^{*}$. Then, by lemma 1 we have $\mathbf{y}_{i}\in X^{\mathbf{v}_{i}^{*}}$ such that
$$\mathbf{v}\in {\rm int}(P_{10})+\mathbf{y}_{i}, \quad i=1,\ 2,\ 5.$$
It is obvious that ${\bf y}_1$, ${\bf y}_2$ and ${\bf y}_5$ are pairwise distinct. Thus, we have
$\varphi(\mathbf{v})\geq3$ and therefore
$$\tau(P_{10})=\varphi(\mathbf{v})+\phi(\mathbf{v})\geq6,\eqno(51)$$
which contradicts the assumption.

\smallskip\noindent
{\bf Subcase 2.2.} {\it $G'= G$.} Then $P_{10}+X^{\mathbf{v}}$ can be divided into two adjacent wheels, as shown by Figure 9.

\begin{figure}[!ht]
\centering
\includegraphics[scale=0.53]{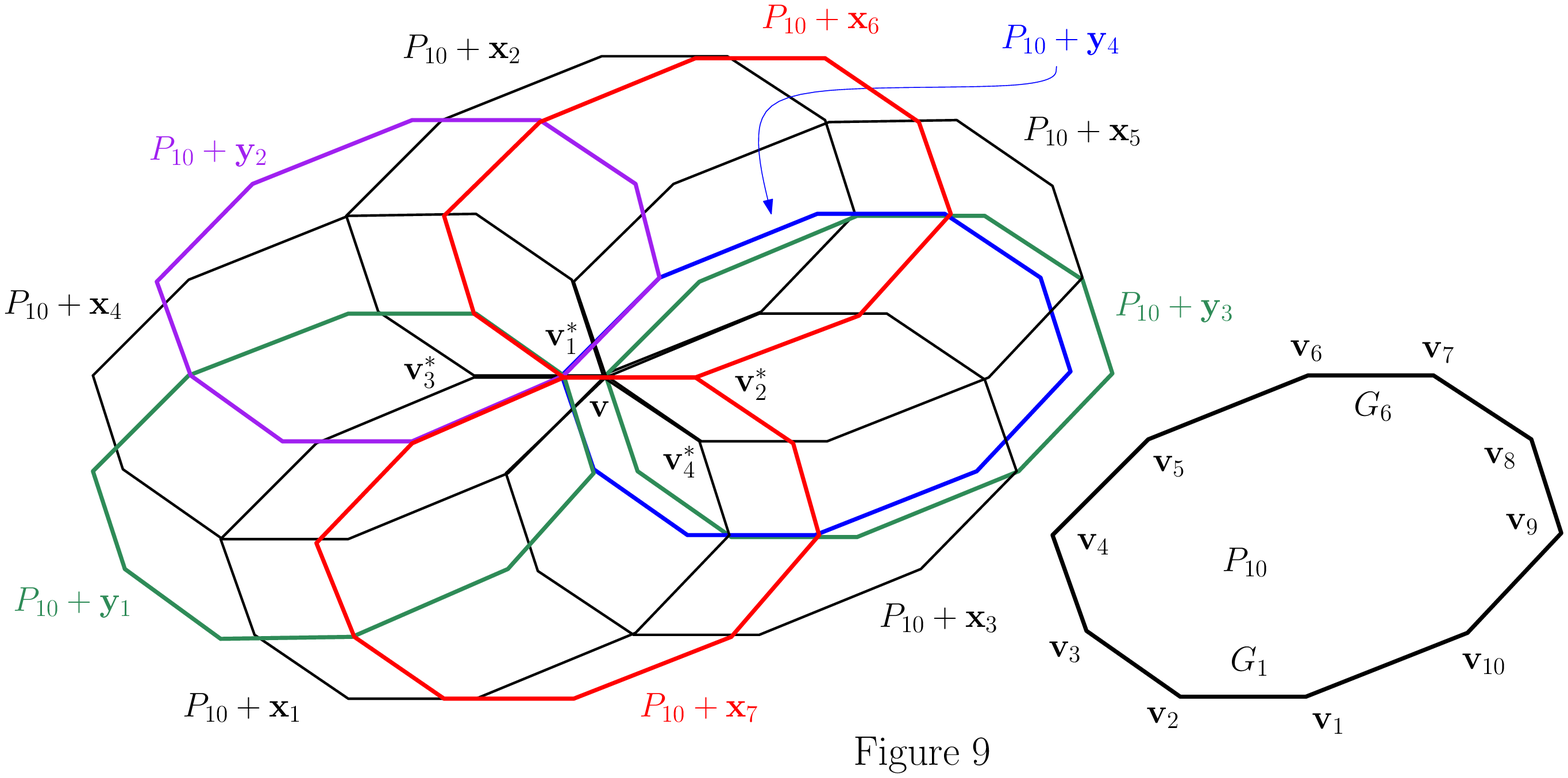}
\end{figure}

Let $P_{10}+\mathbf{x}_{6}$ and $P_{10}+\mathbf{x}_{7}$ be the two translates that contain $G$ as a common edge. Without loss of generality, suppose that $G=G_{6}+\mathbf{x}_{7}$ and $\mathbf{v}=\mathbf{v}_{7}+\mathbf{x}_{1}$, as shown in Figure 9. Let $L$ be the straight line determined by $\mathbf{v}_{1}^{*}$ and $\mathbf{v}_{2}^{*}$, let $G_{1}^{*}$ be the edge of $P_{10}+\mathbf{x}_{1}$ lying on $L$ with ends $\mathbf{v}$ and $\mathbf{v}_{3}^{*}$, and let $G_{2}^{*}$ be the edge of $P_{10}+\mathbf{x}_{1}$ with ends $\mathbf{v}$ and $\mathbf{v}_{4}^{*}$.

It is easy to see that $\phi(\mathbf{v}_{1}^{*})\geq 3$ since ${\bf v}^*_1$ is an interior point of $G^*_1$. In fact, to avoid contradiction, we must have $\phi(\mathbf{v}_{1}^{*})=3$ and the adjacent wheel at $\mathbf{v}_{1}^{*}$ can be divided into two adjacent wheels. Since $\mathbf{v}_{1}^{*}=\mathbf{v}_{6}+\mathbf{x}_{7}$, by Lemma 1 and the structure of the adjacent wheel that consists of five translates, we have three points $\mathbf{y}_{1},$ $\mathbf{y}_{2}$, $\mathbf{y}_{3}\in X^{\mathbf{v}_{1}^{*}}$ such that
$$\mathbf{v}_{1}^{*}=\mathbf{v}_{8}+\mathbf{y}_{1},\quad \mathbf{v}_{3}^{*}\in {\rm int}(P_{10})+\mathbf{y}_1,\eqno(52)$$
$$\mathbf{v}_{1}^{*}=\mathbf{v}_{10}+\mathbf{y}_{2},\quad \mathbf{v}_{3}^{*}\in {\rm int}(P_{10})+\mathbf{y}_{2}\eqno(53)$$
and
$$\mathbf{v}_{1}^{*}=\mathbf{v}_{4}+\mathbf{y}_{3},\quad \mathbf{v}\in {\rm int}(P_{10})+\mathbf{y}_{3}.\eqno(54)$$

Clearly, we also have $\mathbf{v}_{3}^{*}\in {\rm int}(P_{10})+\mathbf{x}_{4}$. Since $\mathbf{v}_{1}^{*}\in {\rm int}(P_{10})+\mathbf{x}_{4}$, thus we have $\mathbf{x}_{4}\notin \{\mathbf{y}_{1},\mathbf{y}_{2}\}$, $\varphi(\mathbf{v}_{3}^{*})\geq3$ and $\phi(\mathbf{v}_{3}^{*})=2$. By Lemma 1 and the structure of the adjacent wheel with five translates, there is a point ${\bf y}_4\in X^{\mathbf{v}_{3}^{*}}$ such that
$$\mathbf{v}_{3}^{*}=\mathbf{v}_{4}+\mathbf{y}_{4},\quad \mathbf{v}\in {\rm int}(P_{10})+\mathbf{y}_{4}.\eqno(55)$$
Furthermore, by Lemma 1 we have $\mathbf{y}_{5}\in X^{\mathbf{v}_{4}^{*}}$ such that $\mathbf{v}\in {\rm int}(P_{10})+\mathbf{y}_{5}.$
By (54), (55) and convexity we have
$$\mathbf{v}_{4}^{*}\in \left({\rm int}(P_{10})+\mathbf{y}_{3}\right)\cap \left({\rm int}(P_{10})+\mathbf{y}_{4}\right),$$ $\mathbf{y}_{5}\notin \{\mathbf{y}_{3},\mathbf{y}_{4}\}$, $\varphi(\mathbf{v})\geq3$ and thus
$$\tau(P_{10})=\varphi(\mathbf{v})+\phi(\mathbf{v})\geq6,\eqno(56)$$
which contradicts the assumption.

\smallskip\noindent
{\bf Subcase 2.3.} $G'\nparallel G$. Suppose that $G$ is an edge of $P_{10}+\mathbf{x}_{6}$ with ends $\mathbf{v}^{*}_{1}$ and $\mathbf{v}^{*}_{2}$ which contains ${\bf v}$ as an interior point. Since $G'\nparallel G$, there is a translates $P_{10}+\mathbf{x}'$ in $X^{\mathbf{v}}$ that meets $P_{10}+\mathbf{x}_{6}$ at a non-singleton part of $G$. Let L be the line determined by  $\mathbf{v}^{*}_{1}$ and $\mathbf{v}^{*}_{2}$. Let $G_{1}^{*}$ be the edge of $P_{10}+\mathbf{x}'$ lying on $L$ with ends $\mathbf{v}^{*}_{3}$ and $\mathbf{v}$. Since $\mathbf{v}^{*}_{1}\in {\rm int}(G_{1}^{*})$, we have $\varphi ({\bf v}^*_1)\ge 2$ and therefore $\phi(\mathbf{v}^{*}_{1})\leq 3$. On the other hand, since $\ell \not= 0$ at ${\bf v}^*_1$, we must have $\phi ({\bf v}^*_1)\ge 3.$ Thus we get
$$\phi ({\bf v}^*_1)=3.\eqno(57)$$

Since $G'\nparallel G$, then $P_{10}+X^{\mathbf{v}}$ cannot be divided into smaller adjacent wheels. There are exact two corresponding inner angles of two translates in $P_{10}+X^{\mathbf{v}}$ at $\mathbf{v}$ that are divided into two positive measure parts by $G^{*}_{1}$. Thus there are exact two translates in $P_{10}+X^{\mathbf{v}}$ that contain both $\mathbf{v}^{*}_{3}$ and $\mathbf{v}^{*}_{1}$ as  interior points. By Lemma 1, there is a translate $P_{10}+\mathbf{y}$ in $P_{10}+X^{\mathbf{v}_{3}^{*}}$ that contains $\mathbf{v}^{*}_{1}$ as an interior point and therefore $\varphi(\mathbf{v}^{*}_{1})\geq3$. Then, by (57) we get
$$\tau(P_{10})=\varphi(\mathbf{v}^*_1)+\phi(\mathbf{v}^*_1)\geq6,\eqno(58)$$
which contradicts the assumption.

As a conclusion, we have proved that $\phi(\mathbf{v})=2$ must hold for all vertices $\mathbf{v}\in V+X$ if $P_{10}+X$ is a five-fold translative tiling.

\medskip
Let $P_{10}$ be a centrally symmetric convex decagon centered at the origin with vertices $\mathbf{v}_{1},$ $\mathbf{v}_{2},$ $\ldots,$ $\mathbf{v}_{10}$ enumerated in the anti-clock order. Let $G_{i}$ denote the edge with ends $\mathbf{v}_{i}$ and $\mathbf{v}_{i+1}$ and let ${\bf u}_i$ denote the middle point of $G_i$. Then, we define
$$\left\{\begin{array}{ll}
{\bf a}_1\hspace{-0.3cm}&={\bf u}_1-{\bf u}_6,\\
{\bf a}_2\hspace{-0.3cm}&={\bf u}_2-{\bf u}_7,\\
{\bf a}_3\hspace{-0.3cm}&={\bf u}_3-{\bf u}_8,\\
{\bf a}_4\hspace{-0.3cm}&={\bf u}_4-{\bf u}_9,\\
{\bf a}_5\hspace{-0.3cm}&={\bf u}_5-{\bf u}_{10}.
\end{array}\right.$$
According to Lemma 8 of Zong \cite{zong}, we have
$${\bf a}_1-{\bf a}_2+{\bf a}_3-{\bf a}_4+{\bf a}_5={\bf o}.\eqno(59)$$

Assume that ${\bf x}_1={\bf o}\in X$. Since $\phi(\mathbf{v})=2$ holds for every vertex $\mathbf{v}\in V+X$, by studying the structure of the adjacent wheel at $\mathbf{v}$ we have
$$\sum_{i=1}^{5} z_{i} \mathbf{a}_{i}\in X, \quad z_{i}\in \mathbb{Z}.$$
For convenience, we define
$$\Lambda =\left\{\sum z_{i} \mathbf{a}_{i}:\  z_{i}\in \mathbb{Z},\ i=1,2,\ldots,5\right\}.\eqno(60)$$

Since $\phi(\mathbf{v})=2$ hold for all vertices, we have $\varphi(\mathbf{v})=3$ for every vertex ${\bf v}$ as well. Suppose that the adjacent wheel at $\mathbf{v}_{1}$ is $P_{10}+\mathbf{x}_{i}$, $i=1,$ $2,$ $\ldots,$ $5$. Let $\mathbf{v}_{i}^{*}$ be the common vertex of $P_{10}+\mathbf{x}_{i}$ and $P_{10}+\mathbf{x}_{i+1}$ other than ${\bf v}_1$ as shown by Figure 10, where $\mathbf{x}_{6}=\mathbf{x}_{1}$ and $\mathbf{x}_{1}={\bf o}$. By Lemma 1, we have $\mathbf{y}_{i}\in X^{\mathbf{v}_{i}^{*}}$ such that $\mathbf{v}_{1}\in {\rm int}(P_{10})+\mathbf{y}_{i}$. In fact, it can be explicitly deduced that
$$\left\{\begin{array}{ll}
\mathbf{y}_{1}\hspace{-0.3cm}&={\bf v}^*_1-{\bf v}_4={\bf a}_{1}-{\bf a}_{4}+{\bf a}_{5},\\
\mathbf{y}_{2}\hspace{-0.3cm}&={\bf v}^*_2-{\bf v}_{10}=2{\bf a}_{1}-{\bf a}_{2}+{\bf a}_{5},\\
\mathbf{y}_{3}\hspace{-0.3cm}&={\bf v}^*_3-{\bf v}_6=2{\bf a}_{1}-2{\bf a}_{2}+{\bf a}_{3},\\
\mathbf{y}_{4}\hspace{-0.3cm}&={\bf v}^*_4-{\bf v}_2=-{\bf a}_{2}+{\bf a}_{3}-{\bf a}_{5},\\
\mathbf{y}_{5}\hspace{-0.3cm}&={\bf v}^*_5-{\bf v}_8=-{\bf a}_{1}+{\bf a}_{2}-{\bf a}_{5}.
\end{array}\right.\eqno(61)
$$

\begin{figure}[!ht]
\centering
\includegraphics[scale=0.51]{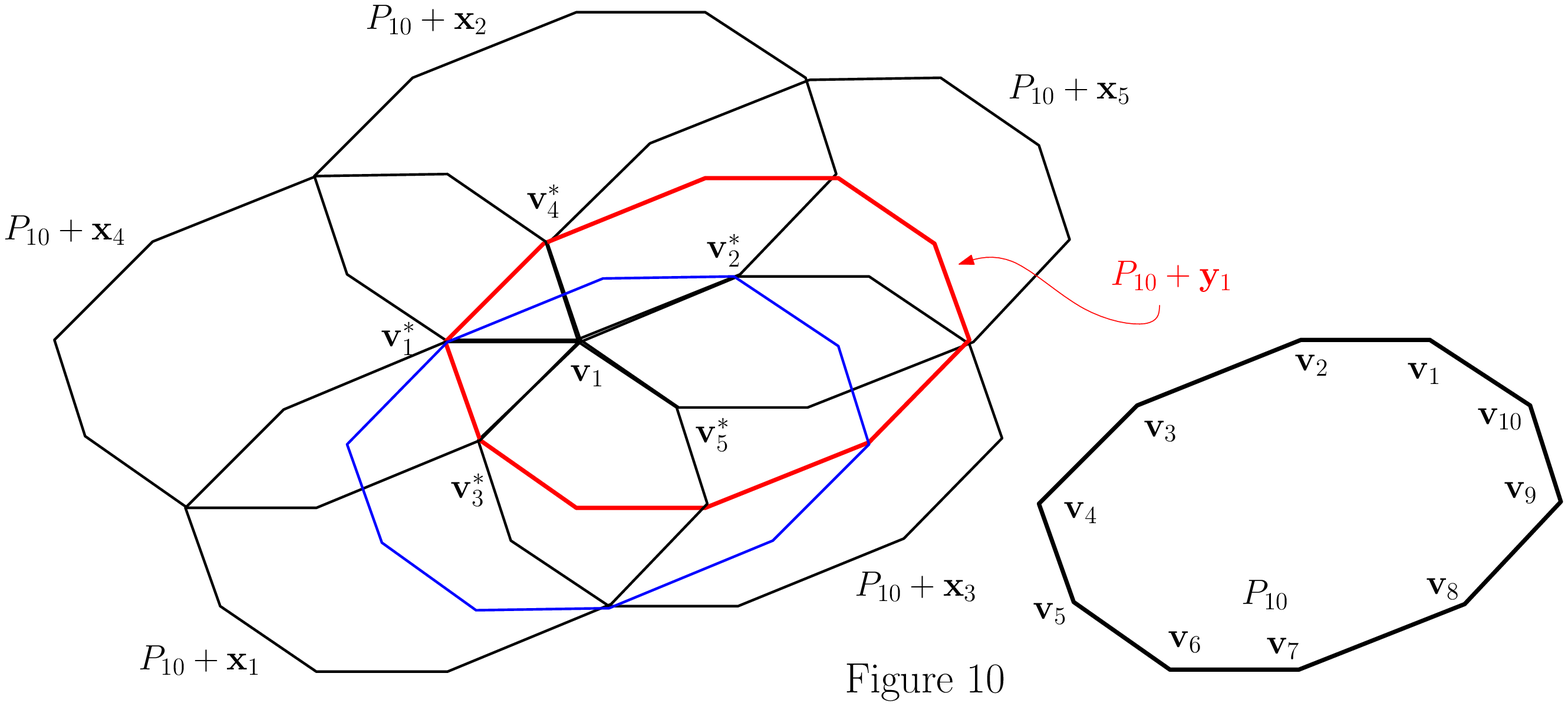}
\end{figure}

By (61) and symmetry it can be shown that $\mathbf{y}_{i}\neq\mathbf{y}_{i+1}$, where $\mathbf{y}_{6}=\mathbf{y}_{1}$. For example, if ${\bf y}_1={\bf y}_2$ (as shown in Figure 10), then by symmetry we will get that $(P_{10}+{\bf x}_2)\cap (P_{10}+{\bf y}_1)$ is a parallelogram and ${\bf y}_1={\bf v}^*_1-{\bf v}_3$, which contradicts the first equation of (61). Thus, any triple of $\{\mathbf{y}_{1},\mathbf{y}_{2},\ldots,\mathbf{y}_{5}\}$ cannot be identical. Restricted by $\varphi(\mathbf{v})=3$, these points have to satisfy one of the following five groups of conditions: {\bf (i).} $\mathbf{y}_{1}=\mathbf{y}_{3}$ and $\mathbf{y}_{2}=\mathbf{y}_{4}$; {\bf (ii).} $\mathbf{y}_{1}=\mathbf{y}_{3}$ and $\mathbf{y}_{2}=\mathbf{y}_{5}$; {\bf (iii).} $\mathbf{y}_{1}=\mathbf{y}_{4}$ and $\mathbf{y}_{2}=\mathbf{y}_{5}$; {\bf (iv).} $\mathbf{y}_{1}=\mathbf{y}_{4}$ and $\mathbf{y}_{3}=\mathbf{y}_{5}$ and {\bf (v).} $\mathbf{y}_{2}=\mathbf{y}_{4}$ and $\mathbf{y}_{3}=\mathbf{y}_{5}$.

\medskip
\noindent{\bf Case (i).} {\it $\mathbf{y}_{1}=\mathbf{y}_{3}$ and $\mathbf{y}_{2}=\mathbf{y}_{4}$}. Then, by (61) and (59) we get
$$\left\{ \begin{array}{ll}
{\bf a}_{1}-{\bf a}_{4}+{\bf a}_{5}=2{\bf a}_{1}-2{\bf a}_{2}+{\bf a}_{3},&\\
2{\bf a}_{1}-{\bf a}_{2}+{\bf a}_{5}=-{\bf a}_{2}+{\bf a}_{3}-{\bf a}_{5},&\\
{\bf a}_1-{\bf a}_2+{\bf a}_3-{\bf a}_4+{\bf a}_5={\bf o},&
\end{array}\right.\eqno(62)$$

$$\left\{ \begin{array}{ll}
2{\bf a}_2-2{\bf a}_{4}+{\bf a}_{5}={\bf a}_{1}+({\bf a}_3-{\bf a}_4),&\\
{\bf a}_4-2{\bf a}_5=2{\bf a}_1-({\bf a}_{3}-{\bf a}_{4}),&\\
{\bf a}_2-{\bf a}_5={\bf a}_1+({\bf a}_3-{\bf a}_4)&
\end{array}\right.$$
and therefore
$$\left\{ \begin{array}{ll}
{\bf a}_1={\bf a}_1,&\\
{\bf a}_2=-2{\bf a}_1+4({\bf a}_3-{\bf a}_4),&\\
{\bf a}_3=-4{\bf a}_1+6({\bf a}_3-{\bf a}_4),&\\
{\bf a}_4=-4{\bf a}_1+5({\bf a}_3-{\bf a}_4),&\\
{\bf a}_5=-3{\bf a}_1+3({\bf a}_3-{\bf a}_4),&
\end{array}\right.\eqno(63)$$
which means that $\Lambda $ is a lattice with a basis $\{ {\bf a}_{1},\ {\bf a}_3-{\bf a}_{4}\}$.
Furthermore, since ${\bf u}_i={1\over 2}{\bf a}_i\in {1\over 2}\Lambda,$ it follows by Lemma 4 that
$P_{10}+\Lambda$ is indeed a five-fold lattice tiling.

\medskip\noindent
{\bf Case (ii).} {\it $\mathbf{y}_{1}=\mathbf{y}_{3}$ and $\mathbf{y}_{2}=\mathbf{y}_{5}$.} Then, by (61) and (59) we have
$$\left\{ \begin{array}{ll}
{\bf a}_{1}-{\bf a}_{4}+{\bf a}_{5}=2{\bf a}_{1}-2{\bf a}_{2}+{\bf a}_{3},&\\
2{\bf a}_{1}-{\bf a}_{2}+{\bf a}_{5}=-{\bf a}_{1}+{\bf a}_{2}-{\bf a}_{5},&\\
{\bf a}_1-{\bf a}_2+{\bf a}_3-{\bf a}_4+{\bf a}_5={\bf o},&
\end{array}\right.\eqno(64)$$

$$\left\{ \begin{array}{ll}
{\bf a}_{1}+{\bf a}_{4}+{\bf a}_{5}=-{\bf a}_3+2({\bf a}_{2}+{\bf a}_5),&\\
3{\bf a}_{1}+4{\bf a}_{5}=2({\bf a}_2+{\bf a}_{5}),&\\
{\bf a}_1-{\bf a}_4+2{\bf a}_5=-{\bf a}_3+({\bf a}_2+{\bf a}_5)&
\end{array}\right.$$
and therefore
$$\left\{ \begin{array}{ll}
{\bf a}_1=8{\bf a}_3-6({\bf a}_2+{\bf a}_5),&\\
{\bf a}_2=6{\bf a}_3-4({\bf a}_2+{\bf a}_5),&\\
{\bf a}_3={\bf a}_3,&\\
{\bf a}_4=-3{\bf a}_3+3({\bf a}_2+{\bf a}_5),&\\
{\bf a}_5=-6{\bf a}_3+5({\bf a}_2+{\bf a}_5),&
\end{array}\right.\eqno(65)$$
which means that $\Lambda $ is a lattice with a basis $\{ {\bf a}_3,\ {\bf a}_2+{\bf a}_5\}$.
Furthermore, since ${\bf u}_i={1\over 2}{\bf a}_i\in {1\over 2}\Lambda,$ it follows by Lemma 4 that
$P_{10}+\Lambda$ is indeed a five-fold lattice tiling.

\medskip
\noindent{\bf Case (iii).} {\it $\mathbf{y}_{1}=\mathbf{y}_{4}$ and $\mathbf{y}_{2}=\mathbf{y}_{5}$.} Then, by (61) and (59) we get
$$\left\{ \begin{array}{ll}
{\bf a}_{1}-{\bf a}_{4}+{\bf a}_{5}=-{\bf a}_{2}+{\bf a}_{3}-{\bf a}_{5},&\\
2{\bf a}_{1}-{\bf a}_{2}+{\bf a}_{5}=-{\bf a}_{1}+{\bf a}_{2}-{\bf a}_{5},&\\
{\bf a}_1-{\bf a}_2+{\bf a}_3-{\bf a}_4+{\bf a}_5={\bf o}&
\end{array}\right.\eqno(66)$$

$$\left\{ \begin{array}{ll}
2{\bf a}_2-{\bf a}_3+2{\bf a}_{5}=-({\bf a}_1-{\bf a}_2)+{\bf a}_4,&\\
{\bf a}_2+2{\bf a}_{5}=-3({\bf a}_{1}-{\bf a}_{2}),&\\
{\bf a}_3+{\bf a}_5=-({\bf a}_1-{\bf a}_2)+{\bf a}_4&
\end{array}\right.$$
and therefore
$$\left\{ \begin{array}{ll}
{\bf a}_1=4{\bf a}_4+6({\bf a}_1-{\bf a}_2),&\\
{\bf a}_2=4{\bf a}_4+5({\bf a}_1-{\bf a}_2),&\\
{\bf a}_3=3{\bf a}_4+3({\bf a}_1-{\bf a}_2),&\\
{\bf a}_4={\bf a}_4,&\\
{\bf a}_5=-2{\bf a}_4-4({\bf a}_1-{\bf a}_2),&
\end{array}\right.\eqno(67)$$
which means that $\Lambda $ is a lattice with a basis $\{ {\bf a}_4,\ {\bf a}_1-{\bf a}_2\}$.
Furthermore, since ${\bf u}_i={1\over 2}{\bf a}_i\in {1\over 2}\Lambda,$ it follows by Lemma 4 that
$P_{10}+\Lambda$ is indeed a five-fold lattice tiling.

\medskip
\noindent
{\bf Case (iv).} {\it $\mathbf{y}_{1}=\mathbf{y}_{4}$ and $\mathbf{y}_{3}=\mathbf{y}_{5}$.} Then, by (61) and (59) we have
$$\left\{ \begin{array}{ll}
{\bf a}_{1}-{\bf a}_{4}+{\bf a}_{5}=-{\bf a}_{2}+{\bf a}_{3}-{\bf a}_{5},&\\
2{\bf a}_{1}-2{\bf a}_{2}+{\bf a}_{3}=-{\bf a}_{1}+{\bf a}_{2}-{\bf a}_{5},&\\
{\bf a}_1-{\bf a}_2+{\bf a}_3-{\bf a}_4+{\bf a}_5={\bf o},&
\end{array}\right.\eqno(68)$$

$$\left\{ \begin{array}{ll}
{\bf a}_2-{\bf a}_{4}+{\bf a}_{5}={\bf a}_3-({\bf a}_1+{\bf a}_{5}),&\\
3{\bf a}_2+2{\bf a}_5={\bf a}_3+3({\bf a}_1+{\bf a}_{5}),&\\
{\bf a}_2+{\bf a}_4={\bf a}_3+({\bf a}_1+{\bf a}_5)&
\end{array}\right.$$
and therefore
$$\left\{ \begin{array}{ll}
{\bf a}_1=4{\bf a}_3-5({\bf a}_1+{\bf a}_5),&\\
{\bf a}_2=3{\bf a}_3-3({\bf a}_1+{\bf a}_5),&\\
{\bf a}_3={\bf a}_3,&\\
{\bf a}_4=-2{\bf a}_3+4({\bf a}_1+{\bf a}_5),&\\
{\bf a}_5=-4{\bf a}_3+6({\bf a}_1+{\bf a}_5),&
\end{array}\right.\eqno(69)$$
which means that $\Lambda $ is a lattice with a basis $\{ {\bf a}_3,\ {\bf a}_1+{\bf a}_5\}$.
Furthermore, since ${\bf u}_i={1\over 2}{\bf a}_i\in {1\over 2}\Lambda,$ it follows by Lemma 4 that
$P_{10}+\Lambda$ is indeed a five-fold lattice tiling.

\medskip
\noindent{\bf Case (v).} {\it $\mathbf{y}_{2}=\mathbf{y}_{4}$ and $\mathbf{y}_{3}=\mathbf{y}_{5}$.} Then, by (61) and (59) we have
$$\left\{ \begin{array}{ll}
2{\bf a}_{1}-{\bf a}_{2}+{\bf a}_{5}=-{\bf a}_{2}+{\bf a}_{3}-{\bf a}_{5},&\\
2{\bf a}_{1}-2{\bf a}_{2}+{\bf a}_{3}=-{\bf a}_{1}+{\bf a}_{2}-{\bf a}_{5},&\\
{\bf a}_1-{\bf a}_2+{\bf a}_3-{\bf a}_4+{\bf a}_5={\bf o},&
\end{array}\right.\eqno(70)$$

$$\left\{ \begin{array}{ll}
2{\bf a}_{1}-{\bf a}_3=-2{\bf a}_{5},&\\
3{\bf a}_{1}+{\bf a}_3-3{\bf a}_4=3({\bf a}_2-{\bf a}_4)-{\bf a}_{5},&\\
{\bf a}_1+{\bf a}_3-2{\bf a}_4=({\bf a}_2-{\bf a}_4)-{\bf a}_5&
\end{array}\right.$$
and therefore
$$\left\{ \begin{array}{ll}
{\bf a}_1=3{\bf a}_5+3({\bf a}_2-{\bf a}_4),&\\
{\bf a}_2=6{\bf a}_5+5({\bf a}_2-{\bf a}_4),&\\
{\bf a}_3=8{\bf a}_5+6({\bf a}_2-{\bf a}_4),&\\
{\bf a}_4=6{\bf a}_5+4({\bf a}_2-{\bf a}_4),&\\
{\bf a}_5={\bf a}_5,&
\end{array}\right.\eqno(71)$$
which means that $\Lambda $ is a lattice with a basis $\{ {\bf a}_5,\ {\bf a}_2-{\bf a}_4\}$.
Furthermore, since ${\bf u}_i={1\over 2}{\bf a}_i\in {1\over 2}\Lambda,$ it follows by Lemma 4 that
$P_{10}+\Lambda$ is indeed a five-fold lattice tiling.

\medskip
As a conclusion of these five cases, Lemma 6 is proved. \hfill{$\Box$}

\medskip\noindent
{\bf Lemma 7 (Zong \cite{zong}).} {\it A convex decagon can form a five-fold lattice tiling of the Euclidean plane if and only if, under a
suitable affine linear transformation, it takes ${\bf u}_1=(0,2)$, ${\bf u}_2=(2,2)$, ${\bf u}_3=(3,1)$, ${\bf u}_4=(3,0)$, ${\bf u}_5=(2,-1)$, ${\bf u}_6=-{\bf u}_1$, ${\bf u}_7=-{\bf u}_2$, ${\bf u}_8=-{\bf u}_3$, ${\bf u}_9=-{\bf u}_4$ and ${\bf u}_{10}=-{\bf u}_5$ as the middle points of its edges.}

\medskip
\noindent
{\bf Remark 2 (Zong \cite{zong}).} Let $W$ denote the quadrilateral with vertices ${\bf w}_1=(-1,2)$, ${\bf w}_2=(-1,{3\over 2})$, ${\bf w}_3=(-{4\over 3}, {4\over 3})$ and ${\bf w}_4=(-{3\over 2},{3\over 2})$.  A centrally symmetric convex decagon with ${\bf u}_1=(0,2)$, ${\bf u}_2=(2,2)$, ${\bf u}_3=(3,1)$, ${\bf u}_4=(3,0)$, ${\bf u}_5=(2,-1)$, ${\bf u}_6=-{\bf u}_1$, ${\bf u}_7=-{\bf u}_2$, ${\bf u}_8=-{\bf u}_3$, ${\bf u}_9=-{\bf u}_4$ and ${\bf u}_{10}=-{\bf u}_5$ as the middle points of its edges if and only if one of its vertices is an interior point of $W$.

\medskip
\noindent
{\bf Proof of Theorem 1.} Theorem 1 follows from Lemma 3, Lemma 5, Lemma 6 and Lemma 7. \hfill{$\Box$}

\vspace{0.6cm}\noindent
{\bf Acknowledgements.} For helpful comments and suggestions, the authors are grateful to Professor S. Robins and Professor G. M. Ziegler. This work is supported by 973 Program 2013CB834201.

\bibliographystyle{amsplain}

\begin{thebibliography}{}
\bibitem{alek}A. D. Aleksandrov, On tiling space by polytopes, {\it Vestnik Leningrad Univ. Ser. Mat. Fiz. Him}. {\bf 9} (1954), 33-43.
\bibitem{boll}U. Bolle, On multiple tiles in $R^2$, {\it Intuitive Geometry,} Colloq. Math. Soc. J. Bolyai {\bf 63}, North-Holland, Amsterdam, 1994.
\bibitem{delo}B. N. Delone, Sur la partition reguli$\grave{e}$re de l'espace $\grave{a}$ $4$ dimensions I, II,
{\it Izv. Akad. Nauk SSSR, Ser. VII} (1929), 79-110; 147-164.
\bibitem{dgsw}M. Dutour Sikiri$\acute{\rm c}$, A. Garber, A. Sch$\ddot{\rm u}$rmann and C. Waldmann, The complete classification of five-dimensional Dirichlet-Voronoi polyhedra of translational lattices,  arXiv:1507.00238.
\bibitem{enge}P. Engel, On the symmetry classification of the four-dimensional parallelohedra, {\it Z. Kristallographie} {\bf 200} (1992), 199-213.
\bibitem{fedo}E. S. Fedorov, Elements of the study of figures, {\it Zap. Mineral. Imper. S. Petersburgskogo Ob$\check{s}$$\check{c}$},
{\bf 21}(2) (1885), 1-279.
\bibitem{furt}P. Furtw\"angler, \"Uber Gitter konstanter Dichte, {\it Monatsh. Math. Phys.} {\bf 43} (1936), 281-288.
\bibitem{grs}N. Gravin, S. Robins and D. Shiryaev, Translational tilings by a polytope, with multiplicity. {\it Combinatorica}
{\bf 32} (2012), 629-649.
\bibitem{gkrs}N. Gravin, M. N. Kolountzakis, S. Robins and D. Shiryaev, Structure results for multiple tilings in 3D.
{\it Discrete Comput. Geom.} {\bf 50} (2013), 1033-1050.
\bibitem{grub}P. M. Gruber and C. G. Lekkerkerker, {\it Geometry of Numbers} (2nd ed.), North-Holland, 1987.
\bibitem{hajo}G. Haj\'os, \"Uber einfache und mehrfache Bedeckung des $n$-dimensionalen Raumes mit einem W\"urfelgitter, {\it Math. Z.}
{\bf 47} (1941), 427-467.
\bibitem{hilb}D. Hilbert, Mathematical Problems, {\it G\"ottinger Nachr.} (1900), 253-297. {\it Bull. Amer. Math. Soc.} {\bf 8} (1902), 437-479;
{\bf 37} (2000), 407-436.
\bibitem{kolo}M. N. Kolountzakis, On the structure of multiple translational tilings by polygonal regions. {\it Discrete Comput. Geom}.
{\bf 23} (2000), 537-553.
\bibitem{lazo}J. C. Lagarias and C. Zong, Mysteries in packing regular tetrahedra, {\it Notices Amer. Math. Soc.} {\bf 59} (2012), 1540-1549.
\bibitem{mann}C. Mann, J. McLoud-Mann and D. Von Derau, Convex pentagons that admit $i$-block transitive tilings, {\it Geom. Dedicata} (2018), in press.
\bibitem{mcmu}P. McMullen, Convex bodies which tiles space by translation, {\it Mathematika} {\bf 27} (1980), 113-121.
\bibitem{mink}H. Minkowski, Allgemeine Lehrs\"atze \"uber konvexen Polyeder, {\it Nachr. K. Ges. Wiss. G\"ottingen, Math.-Phys. KL}. (1897), 198-219.
\bibitem{rao}M. Rao, Exhaustive search of convex pentagons which tile the plane, arXiv:1708.00274
\bibitem{rein}K. Reinhardt, \"Uber die Zerlegung der Ebene in Polygone, {\it Dissertation,} Universit\"at Frankfurt am Main, 1918.
\bibitem{robi}R. M. Robinson, Multiple tilings of $n$-dimensional space by unit cubes, {\it Math. Z.} {\bf 166} (1979), 225-275.
\bibitem{stog}M. I. $\check{S}$togrin, Regular Dirichlet-Voronoi partitions for the second triclinic group (in Russian),
{\it Proc. Steklov. Inst. Math.} {\bf 123} (1975).
\bibitem{venk}B. A. Venkov, On a class of Euclidean polytopes, {\it Vestnik Leningrad Univ. Ser. Mat. Fiz. Him}. {\bf 9} (1954), 11-31.
\bibitem{voro}G. F. Voronoi,  Nouvelles applications des paramm\`etres continus \`a la th\'eorie des formes quadratiques.
Deuxi\`eme M\'emoire. Recherches sur les parall\'elo\`edres primitifs, {\it J. reine angew. Math.} {\bf 134} (1908), 198-287; {\bf 135} (1909), 67-181.
\bibitem{yz1}Q. Yang and C. Zong, Multiple lattice tiling in Euclidean spaces, {\it Bull. London Math. Soc.} in press.
\bibitem{yz2}Q. Yang and C. Zong, Multiple translative tiling in Euclidean spaces, arXiv:1711.02514.
\bibitem{zong05}C. Zong, What is known about unit cubes. {\it Bull. Amer. Math. Soc.} {\bf 42} (2005), 181-211.
\bibitem{zong06}C. Zong, {\it The Cube: A Window to Convex and Discrete Geometry.} Cambridge University Press, Cambridge, 2006.
\bibitem{zong14}C. Zong, Packing, covering and tiling in two-dimensional spaces, {\it Expo. Math.} {\bf 32} (2014), 297-364.
\bibitem{zong}C. Zong, Characterization of the two-dimensional five-fold lattice tiles, arXiv:1712.01122.

\end{thebibliography}
{}

\vspace{0.8cm}
\noindent
Qi Yang, School of Mathematical Sciences, Peking University, Beijing 100871, China

\medskip
\noindent
Corresponding author:

\noindent
Chuanming Zong, Center for Applied Mathematics, Tianjin University, Tianjin 300072, China

\noindent
Email: cmzong@math.pku.edu.cn

\end{document}